\documentclass[10pt,a4paper,leqno]{article}
\usepackage{amsmath,amssymb,amscd,latexsym}

\usepackage[all]{xy}

\newcommand{\A}{{\mathbb{A}}}
\newcommand{\C}{{\mathbb{C}}}
\newcommand{\F}{{\mathbb{F}}}
\newcommand{\Ge}{\mathbb{G}}
\newcommand{\Lee}{\mathbb{L}}
\newcommand{\N}{{\mathbb{N}}}
\newcommand{\Q}{{\mathbb{Q}}}

\newcommand{\Pro}{{\mathbb{P}}}

\newcommand{\Z}{{\mathbb{Z}}}
\newcommand{\Zl}{{\mathbb{Z}}/\ell}
\newcommand{\Zln}{{\mathbb{Z}}/\ell^{\nu}}

\newcommand{\car}{\mathrm{char}\;}
\newcommand{\ddeg}{\mathrm{deg}}
\newcommand{\et}{\mathrm{\acute{e}t}}

\newcommand{\Ev}{\mathrm{Ev}}
\newcommand{\Ho}{\mathrm{Ho}}
\newcommand{\id}{\mathrm{id}}
\newcommand{\ok}{\overline{k}}
\newcommand{\Xok}{X_{\ok}}

\newcommand{\hOmega}{\hat{\Omega}}
\newcommand{\pro}{\mathrm{pro}}

\newcommand{\Spec}{\mathrm{Spec}\,}
\newcommand{\colim}{\mathrm{colim}}

\newcommand{\Et}{\mathrm{Et}\,}
\newcommand{\hEt}{\hat{\mathrm{Et}}\,}

\newcommand{\hEtsp}{\hat{\mathrm{Et}}_\mathrm{Sp}\,}
\newcommand{\thEtsp}{\tilde{\hat{\mathrm{Et}}}_\mathrm{Sp}\,}
\newcommand{\LhEt}{\mathrm{L}\hat{\mathrm{Et}}\,}
\newcommand{\LthEt}{\mathrm{L}\tilde{\hat{\mathrm{Et}}}\,}
\newcommand{\LhEtsp}{\mathrm{L}\hat{\mathrm{Et}}_\mathrm{Sp}\,}

\newcommand{\compl}{\hat{(\cdot)}}

\newcommand{\Gal}{\mathrm{Gal}}

\newcommand{\Hom}{\mathrm{Hom}}
\newcommand{\homp}{\mathrm{hom}_{\ast}}

\newcommand{\Map}{\mathrm{Map}}

\newcommand{\sk}{\mathrm{sk}}

\newcommand{\Sm}{\mathrm{Sm}}
\newcommand{\Thom}{\mathrm{Th}}

\newcommand{\thomet}{\mathrm{th}_{\et}}

\newcommand{\Ch}{{\mathcal C}}
\newcommand{\Chk}{{\mathcal C}/\hEt k}

\newcommand{\Eh}{{\mathcal E}}
\newcommand{\hEh}{\hat{\mathcal E}}
\newcommand{\Fh}{{\mathcal F}}

\newcommand{\Hh}{{\mathcal H}}

\newcommand{\hHh}{\hat{{\mathcal H}}}
\newcommand{\hHhp}{\hat{{\mathcal H}}_{\ast}}
\newcommand{\Kh}{\mathcal{K}}
\newcommand{\LKC}{L_{\Kh}\mathcal{C}}
\newcommand{\LeU}{L_{\et}U(k)}

\newcommand{\LU}{LU(k)}

\newcommand{\MVh}{\mathcal{MV}_k}

\newcommand{\Oh}{{\mathcal O}}

\newcommand{\Rh}{{\mathcal R}}
\newcommand{\Sh}{{\mathcal S}}

\newcommand{\SHh}{{\mathcal{SH}}}
\newcommand{\SHhS}{{\mathcal{SH}^{S^1}(k)}}
\newcommand{\hSHh}{{\hat{\mathcal{SH}}}}

\newcommand{\hSHhk}{{\hat{\mathcal{SH}}}/\hEt k}
\newcommand{\hSHhtk}{{\hat{\mathcal{SH}}}_2/\hEt k}
\newcommand{\hSh}{\hat{\mathcal S}}
\newcommand{\hShk}{\hat{\mathcal S}/\hEt k}

\newcommand{\hShp}{\hat{\mathcal S}_{\ast}}
\newcommand{\hShpk}{\hat{\mathcal S}_{\ast}/\hEt k}
\newcommand{\hSp}{\mathrm{Sp}(\hShp)}
\newcommand{\hSpk}{\mathrm{Sp}(\hShp/\hEt k)}
\newcommand{\Sp}{\mathrm{Sp}(\Sh_{\ast})}
\newcommand{\SpC}{\mathrm{Sp}(\mathcal{C},T)}

\newcommand{\SpMVP}{\mathrm{Sp}(\MVh,\Pro^1 \wedge \cdot)}

\newcommand{\SpLUQ}{\mathrm{Sp}(\LU,Q \wedge \cdot)}

\newcommand{\Xh}{\mathcal{X}}

\newcommand{\hKU}{\hat{KU}}

\newcommand{\hMU}{\hat{MU}}

\newcommand{\tP}{\tilde{P}}
\newcommand{\tQ}{\tilde{Q}}

\newtheorem{theorem}{Theorem}[section]
\newtheorem{lemma}[theorem]{Lemma}
\newtheorem{prop}[theorem]{Proposition}
\newtheorem{defn}[theorem]{Definition}

\newtheorem{cor}[theorem]{Corollary}

\newtheorem{example}[theorem]{Example}

\newtheorem{remark}[theorem]{Remark}
\newtheorem{condtheorem}[theorem]{Conditional Theorem}

\newenvironment{proof}{\noindent {\bf Proof}}{\mbox{}\hspace*{\fill}$\Box$}
\begin{document}
\title{Stable \'etale realization and \'etale cobordism}
\author{Gereon Quick}
\date{}
\maketitle
\begin{abstract}
We show that there is a stable homotopy theory of profinite spaces and use it for two main applications. On the one hand we construct an \'etale topological realization of the stable $\A^1$-homotopy theory of smooth schemes over a base field of arbitrary characteristic in analogy to the complex realization functor for fields of characteristic zero.\\ 
On the other hand we get a natural setting for \'etale cohomology theories. In particular, we define and discuss an \'etale topological cobordism theory for schemes. It is equipped with an Atiyah-Hirzebruch spectral sequence starting from \'etale cohomology. Finally, we construct maps from algebraic to \'etale cobordism and discuss algebraic cobordism with finite coefficients over an algebraically closed field after inverting a Bott element.
\end{abstract}

\section{Introduction}
For the proof of the Milnor Conjecture Voevodsky invented algebraic cobordism, a new cohomology theory for schemes represented by the Thom spectrum in his new framework of $\A^1$-homotopy theory. Using the known topological realization functor for schemes over a field of characteristic zero he linked his theory to the classical one by showing that the complex topological realization of the algebraic cobordism theory yields a map to complex cobordism. \\
Independently, Levine and Morel constructed in a geometric way an algebraic cobordism theory as the universal oriented cohomology theory on smooth schemes and used it to prove Rost's Degree Formula in characteristic zero.\\ 
Looking at the following diagram of cohomology theories
$$\begin{array}{ccc}
\mbox{algebraic cobordism} & \leftrightarrow  & ?\\
\updownarrow &   & \updownarrow\\
\mbox{algebraic K-Theory} & \leftrightarrow & \mbox{\'etale K-Theory}\\
\updownarrow &   & \updownarrow\\
\mbox{motivic cohomology} & \leftrightarrow & \mbox{\'etale cohomology}
\end{array}$$
gives rise to the following questions. Is there an \'etale cobordism theory that fits into this picture? Can the \'etale theories be described in some kind of stable homotopy theory? Can one replace the complex realization functor in characteristic zero by a topological realization functor on the stable motivic category for arbitrary characteristics? Is there a statement similar to the result of Thomason that algebraic and \'etale cobordism with finite coefficients agree after inverting a Bott element? In this paper we give an answer to these questions.\\
Sending a scheme to its \'etale topological type, first constructed by Artin/ Mazur and Friedlander, gives naturally rise to cohomology theories. The idea is due to Eric Friedlander who constructed in \cite{etaleK} a first version of an \'etale topological K-theory for schemes, which had been developed further by Dwyer and Friedlander in \cite{dwyfried}. This theory turned out to be a powerful tool for the study of algebraic K-theory with finite coefficients. In particular, Thomason proved in his famous paper \cite{thomason} that algebraic K-theory with finite coefficients agrees with \'etale K-theory after inverting a Bott element.\\ 
\v{C}ech and \'etale homotopy theory of schemes has also been studied by Cox and Edwards/Hastings in the 1970s. At the end of the 1970s, in \cite{snaith} Victor P. Snaith had constructed a $p$-adic cobordism theory for schemes. His approach is close to the definition of algebraic K-theory by Quillen. He defines for every scheme $V$ over $\F_q$, a topological cobordism spectrum $\underline{A\F}_{q,V}$. The homotopy groups of this spectrum are the $p$-adic cobordism groups of $V$. He has calculated these groups for projective bundles, Severi-Brauer schemes and other examples. \\
In the beginning of the 1980s, Roy Joshua had already studied a version of \'etale (co-)bordism as the generalized (co-)homology theory represented by a mod $\ell$-variant of $MU$ on the \'etale homotopy type as in Friedlander's book \cite{fried}, where Friedlander also shows the existence of an Atiyah-Hirzebruch spectral sequence for such theories. He had also defined an algebraic bordism theory in analogy to the topological bordism relation and constructed a map from algebraic to \'etale bordism using the Thom-Pontrjagin construction and tubular neighbourhoods for geometrically unibranched projective varietes over an algebraically closed field of arbitrary characteristic. Joshua did not publish these ideas, since the interest in algebraic cobordism only arose about ten years later. But his geometric constructions and his application to a Spanier-Whitehead duality in \'etale homotopy theory can be found in \cite{joshua}.\\      
In the present paper, we also follow Friedlander's idea, but in a slightly different setting. We consider the profinite completion $\hEt$ of Friedlander's \'etale topological type functor of \cite{fried} with values in the category of simplicial profinite sets $\hSh$. Based on the work of Hovey \cite{hovey} and Morel \cite{ensprofin}, we construct a stable homotopy category $\hSHh$ for simplicial profinite sets for any fixed prime number $\ell$. This is one possible category for the study of \'etale topological cohomology theories. We define general cohomology theories such as \'etale K-Theory, \'etale Morava K-theory and \'etale cobordism as the theories represented by profinite spectra such as $\hat{KU}$, $\hat{K(n)}$ and $\hMU$ respectively. This approach unifies different \'etale theories and coincides with the well known theories of \'etale cohomology and \'etale K-Theory with finite coefficients over a separably closed field, while \'etale cobordism is a new theory, which will  be discussed in this paper.\\ 
An advantage of this approach is the immediate construction of an Atiyah-Hirzebruch spectral sequence for smooth schemes over a field starting from \'etale cohomology and converging to \'etale cobordism with finite coefficients.\\
For a base field without $\ell$-torsion and with finite $\ell$-cohomological dimension, for example a separably closed, a finite or a local field, \'etale cobordism in even degrees is an oriented cohomology theory, i.e. it has Chern classes. Since the algebraic cobordism of Levine and Morel \cite{lm} is universal for such functors, there is a canonical map from this algebraic to \'etale cobordism in even degrees.\\
For the comparison with Voevodsky's theory \cite{a1hom}, we have to extend the \'etale topological type functor to the stable homotopy category of smooth schemes over an arbitrary base field. Isaksen \cite{a1real} and Schmidt \cite{schmidt} extended independently the \'etale topological type functor to the unstable $\A^1$-homotopy category. One of our main results is the existence of an \'etale realization on the {\em stable} $\A^1$-homotopy category:
\begin{theorem}
Let $k$ be an arbitrary base field. The functor $\hEt$ induces an \'etale realization of the stable motivic homotopy category of $\Pro^1$-spectra: $$\LhEt: \SHh(k) \to \hSHhtk.$$
\end{theorem}
In the theorem $\hSHhtk$ denotes the stable homotopy category of profinite $S^2$-spectra over $\hEt k$. In particular, the functor $\LhEt$ is equal to $\hEt$ on the suspension spectra of smooth schemes and sends stable $\A^1$-equivalences between suspension spectra of smooth schemes to isomorphisms in $\hSHhtk$. Over a separably closed base field the image of the motivic Thom spectrum is even canonically isomorphic in the stable profinite homotopy category to the profinite $\hat{MU}$-spectrum.\\
If $k$ has no $\ell$-torsion and has finite $\ell$-cohomological dimension, this realization yields a canonical map from algebraic to \'etale cobordism. The above map fits into a commutative diagram with the map from the Levine-Morel theory.\\
Furthermore, we show that the absolute Galois group of the base field acts trivially on the cobordism of the separable closure of the field. Together with the Atiyah-Hirzebruch spectral sequence this allows us to determine the \'etale cobordism of local fields. \\
Finally, there is a candidate for a Bott element in Voevodsky's theory that is mapped to an invertible element in \'etale cobordism. Together with the Atiyah-Hirzebruch spectral sequence for \'etale cobordism and the announced motivic Atiyah-Hirzebruch spectral sequence for algebraic cobordism of Hopkins and Morel, the results of Levine \cite{bott} can be used to show that, at least over an algebraically closed field, algebraic and \'etale cobordism with finite coefficients agree after inverting a Bott element:
\begin{condtheorem}\label{Bottconjecture2}
Let $X$ be a smooth scheme of finite type over an algebraically closed field $k$ with $\car k \neq \ell$. Let $\beta \in MGL^{0,1}(k;\Zln)$ be the Bott element. If we assume the existence and convergence of an Atiyah-Hirzebruch spectral sequence from motivic cohomology to algebraic cobordism, then there is an isomorphism
$$\phi: (\oplus_{p,q}MGL^{p,q}(X;\Zln))[\beta ^{-1}] \stackrel{\cong}{\longrightarrow} \oplus_{p,q}\hMU_{\et}^{p}(X;\Zln).$$
\end{condtheorem}
A similar result should hold over an arbitrary base field. But over a base field which is not separably closed, the proposed definition of \'etale cobordism might not be the best possible. The problem is that this homotopical construction does not reflect the twists in \'etale cohomology. In order to solve this problem, Morel suggested to the author to study maps in the homotopy theory relative over the base $\hEt k$ and to define \'etale cobordism via $\hEt(MGL)$ instead of $\hMU$. Over a separably closed field both approaches agree. The author hopes to be able to discuss this in future work.\\ 
The content of the paper is as follows. In the first section we construct the stable homotopy category $\hSHh$ of profinite spaces and discuss generalized cohomology theories of profinite spaces. In the next section, $\hSHh$ will turn out to be the target category for the stable \'etale topological realization functor. It also yields a natural frame work for the study of generalized \'etale cohomology theories such as \'etale cobordism, which we discuss in the fourth section. In the last part we will compare algebraic and \'etale cobordism and will discuss the conditional theorem above.\\
I would like to use this opportunity to thank first of all Christopher Deninger for the suggestion to work on this topic. He helped me with important suggestions and gave me the motivating impression to work on something that captured his interest. I am grateful to Marc Levine, Fabien Morel and Alexander Schmidt, who answered all my questions and explained to me both the great picture and the little details. I thank Francois-Xavier Dehon for an answer to a specific problem on K\"unneth isomorphisms for cobordism with coefficients. Finally, I would like to thank the kind referee who made several detailed and helpful suggestions.

\section{Profinite spectra}
\subsection{Profinite spaces}
We want to construct the target category for the stable \'etale realization, the stable homotopy category of profinite spaces. Therefore, we recall some facts about profinite spaces and their homotopy groups.\\ 
Let $\Eh$ denote the category of sets and let $\Fh$ be the full subcategory of $\Eh$
whose objects are finite sets. Let $\hEh$ be the category of compact and totally disconnected
topological spaces. We may identify $\Fh$ with a full subcategory of $\hEh$ in the obvious way. The limit functor $\lim :\pro-\Fh \to \hEh$, which sends a pro-object $X$ of $\Fh$ to the limit in $\hEh$ of the diagram corresponding to $X$, is an equivalence of categories.\\
We denote by $\hSh$ (resp. $\Sh$) the category of simplicial profinite sets (resp. simplicial sets). The objects of $\hSh$ (resp. $\Sh$) will be called {\em profinite spaces} (resp. {\em spaces}). The forgetful functor $\hEh \to \Eh$ admits a left adjoint $\compl:\Eh \to \hEh$. It induces a functor $\compl:\Sh \to \hSh$, which is called {\em profinite completion}. It is left adjoint to the forgetful functor $|\cdot|:\hSh \to \Sh$ which sends a profinite space to its underlying simplicial set.\\ 
For a profinite space $X$ we define the set $\Rh(X)$ of simplicial open equivalence relations on $X$. An element $R$ of $\Rh(X)$ is a simplicial profinite subset of the product $X\times X$ such that, in each degree $n$, $R_n$ is an equivalence relation on $X_n$ and an open subset of $X_n\times X_n$. It is ordered by inclusion. For every element $R$ of $\Rh(X)$, the quotient $X/R$ is a simplicial finite set and the map $X \to X/R$ is a map of profinite spaces. When we consider the limit $\lim_{R\in \Rh(X)} X/R$ in $\hSh$, the map $X \to \lim_{R\in \Rh(X)} X/R$ is an isomorphism, cf. \cite{ensprofin}, Lemme 1.\\
Let $X$ be a profinite space. The continuous cohomology $H^{\ast}(X;\pi)$ of $X$ with coefficients in the profinite abelian group $\pi$ is defined as the cohomology of the complex $C^{\ast}(X;\pi)$ of continuous cochains of $X$ with values in $\pi$, i.e. $C^n(X;\pi)$ denotes the set $\Hom_{\hat{\Eh}}(X_n,\pi)$ of continuous maps $\alpha:X_n \to \pi$ and the differentials $\delta^n:C^n(X;\pi)\to C^{n+1}(X;\pi)$ are the morphisms associating to $\alpha$ the map $\sum_{i=0}^{n+1}\alpha \circ d_i$, where $d_i$ denotes the $i$th face map of $X$. If $\pi$ is a finite abelian group and $Z$ a simplicial set, then the cohomologies $H^{\ast}(Z;\pi)$ and $H^{\ast}(\hat{Z};\pi)$ are canonically isomorphic. \\
Let $\ell$ be a fixed prime number. Fabien Morel has shown in \cite{ensprofin} that the category $\hSh$ can be given the structure of a closed model category. The weak equivalences are the maps inducing isomorphisms in continuous cohomology with coefficients $\Zl$; the cofibrations are the degreewise monomorphisms. The homotopy category is denoted by $\hHh$. Morel has also given an explicit construction of fibrant replacements in $\hSh$, cf. \cite{ensprofin}, 2.1. It is based on the $\Zl$-completion functor of \cite{bouskan}. We denote the fibrant replacement in $\hSh$ of a profinite space $X$ by $\hat{X}^{\ell}$ and call it the $\ell$-completion of $X$. If $Y$ is a simplicial set, we also denote by $\hat{Y}^{\ell}$ the $\ell$-completion of its profinite completion $\hat{Y}$.\\ 
The category $\hSh$ has a pointed analogue $\hShp$ whose objects are maps $\ast \to X$ in $\hSh$, where $\ast$ denotes the constant simplicial set equal to a point. Its morphisms are maps in $\hSh$ that respect the basepoints. The forgetful functor $\hShp \to \hSh$ has a left adjoint, which consists in adding a disjoint basepoint $X \mapsto X_+$. The category $\hShp$ has the obvious induced closed model category structure. The product for two pointed profinite spaces $X$ and $Y$ in $\hShp$ is the smash product $X \wedge Y \in \hShp$, defined in the usual way as the quotient $(X \times Y)/(X \vee Y)$ in $\hShp$. One can show that (smash) products are compatible with $\ell$-completion in $\hHh$ (resp. $\hHhp$).\\ 
Morel has shown that the categories $\hSh$ and $\hShp$ have natural simplicial structures, see \cite{ensprofin}. In particular, for a pointed profinite space $Y$, there are profinite mapping spaces $\homp(W,Y)$ for every pointed simplicial finite set $W$. The simplicial profinite set $\homp(W,Y)$ is defined in degree $n$ by the set $\Hom_{\hSh}(W\wedge \Delta[n]_+,Y)$, whose profinite structure is inherited from $Y$ by considering it as the limit over the integers $s$ and the elements $R$ of $\Rh(Y)$ of the finite sets $\Hom_{\hSh}(\sk_s W\wedge \Delta[n]_+,Y/R)$. Here, $\sk_sW$ denotes the $s$-skeleton of $W$. As an application, let $S^1$ be the simplicial finite set $\Delta[1]/\partial \Delta[1]$. As usual we set $\hat{\Omega}X:=\homp(S^1,X)$. In particular, $\hat{\Omega}X$ is fibrant if $X$ is fibrant. This adjunction may be extended to the homotopy category $\hHhp$. For every profinite space $X$ and every profinite space $Y$, there is a natural bijection $\Hom_{\hHhp}(S^1 \wedge X, Y) \cong \Hom_{\hHhp}(X,\hat{\Omega}Y))$. Since $S^1\wedge X$ is a cogroup object and $\hOmega X$ a group object in $\hHhp$, we conclude that the previous bijection of maps in $\hHhp$ is in fact an isomorphism of groups, cf. \cite{dehon}. \\
There are several possibilities to define the homotopy groups of a pointed profinite space. One could define $\pi_kX$ to be the homotopy group $\pi_k|\hat{X}^{\ell}|$ of the underlying simplicial set of the fibrant replacement in $\hShp$ of $X$. Dehon suggests the following approach in \cite{dehon} Section 1.3. If $X$ is a pointed fibrant profinite space, we define $\pi_0X$ to be the coequalizer in $\hEh$ of the diagram $d_0,d_1:X_1 \stackrel{\to}{\to} X_0$;  and $\pi_kX$ for $k\geq 1$ is defined to be the group $\pi_0\hat{\Omega}^kX$. One can easily see that $\pi_k X$ has a natural structure of a pro-$\ell$-group for every $k\geq 1$. Both definitions agree since $\pi_k(\hat{\Omega}X)=\pi_{k+1}X$.\\
We conclude this discussion with a collection of results on the homotopy groups of the $\ell$-completion $\hat{X}^{\ell} \in \hSh$ of a simplicial set $X$, which can be deduced from the methods of \cite{bouskan} VI, \S 5. 
\begin{prop}\label{piofZl-completion}
Let $X$ be a pointed connected simplicial set.\\ 
1. We have an isomorphism 
$$\pi_1(\hat{X}^{\ell}) \cong \hat{\pi_1(X)}^{\ell}.$$
In particular, if $\pi_1X$ is a finitely generated abelian group, there is an isomorphism $\pi_1(\hat{X}^{\ell}) \cong  \Z_{\ell} \otimes_{\Z} \pi_1X$, where $\Z_{\ell}$ denotes the $\ell$-adic integers.\\
2. If $X$ is in addition simply connected, and its higher homotopy groups are finitely generated, then we have for all $n \geq 2$
$$\pi_n \hat{X}^{\ell} \cong \hat{\pi_n X}^{\ell} \cong \Z_{\ell} \otimes_{\Z} \pi_nX$$
\end{prop}

\paragraph{The model structure on $\hSh$ is fibrantly generated}
The main technical result of this subsection is that the model structure on $\hSh$, hence also on $\hShp$, is fibrantly generated. We recall some notations and constructions from \cite{ensprofin}. Let $n \geq 0$ be a non-negative integer and $S$ be a profinite set. The functor $\hSh^{\mathrm{op}}\to \Eh$, $X \mapsto \Hom_{\hEh}(X_n,S)$ is representable by a simplicial profinite set, which is denoted by $L(S,n)$. It is given by the formula
$$L(S,n):\Delta^{\mathrm{op}} \to \hat{\Eh} , ~[k] \mapsto S^{\Hom_{\Delta}([n],[k])}.$$
If $M$ is a profinite abelian group, then $L(M,n)$ has a natural structure of a simplicial profinite abelian group and the abelian group $\Hom_{\hSh}(X,L(M,n))$ can be identified with the group $C^n(X;M)$ of continuous $n$-cochains with values in $M$. Furthermore, for every $k$, $L(M,\ast)([k])$ may be considered in the usual way as an abelian cochain complex. \\
For a profinite abelian group $M$, let $Z^n(X;M)$ denote the abelian group of $n$-cocycles of the complex $C^{\ast}(X;M)$. The functor $\hSh^{\mathrm{op}} \to \Eh$, $X \mapsto Z^n(X;M)$ is also representable by a simplicial profinite abelian group, which is denoted by $K(M,n)$, called the profinite Eilenberg-Mac Lane space of type $(M,n)$. The homomorphism $C^n(X;M)\to Z^{n+1}(X;M)$ given by the differential defines a natural map of simplicial profinite abelian groups $L(M,n) \to K(M,n+1)$. 
Consider the two sets of morphisms of Lemme 2 in \cite{ensprofin}:
$$P:=\{L(M,n)\to K(M,n+1), ~ K(M,n)\to \ast | M ~ \mbox{is an abelian pro-$\ell$-group}, n\geq 0 \}$$ 
and
$$Q:=\{L(M,n)\to \ast |  M~\mbox{is an abelian pro-$\ell$-group}, n\geq 0\}.$$
\begin{theorem}\label{fibgen}
The simplicial model structure on $\hSh$, in which the weak equivalences are the 
$\Z/\ell$-cohomological isomorphisms and the cofibrations are the dimensionwise monomorphisms, is left proper and fibrantly generated with $P$ as the set of generating fibrations and $Q$ as the set of generating trivial fibrations.
\end{theorem}
\begin{proof}
The left properness is due to the fact that all objects in $\hSh$ are cofibrant, see for example \cite{hirsch} Corollary 13.1.3. The assertion of the theorem now follows from the methods of the proofs of Lemme 2 and Proposition 3 in \cite{ensprofin}. One can check for example the list of dual criteria for a cofibrantly generated model structure of Kan's Theorem 11.3.1 in \cite{hirsch}, where we only need the existence of small limits since we use cosmallness. One can find a detailed proof in \cite{pec}.
\end{proof}
\subsection{Profinite spectra}
\paragraph{The stable model structure}
The well known stabilization of the category of simplicial spectra of \cite{bousfried}
uses the fact that $\Sh$ is proper in an essential way in order to construct functorial factorizations. 
Since the category $\hSh$ is only left but not right proper, we have to find a different method. Hovey \cite{hovey} has pointed out a general way to stabilize a left proper cellular model category with respect to a left Quillen endofunctor $T$.\\
We have to start with a result on the localization of model categories. We recall the definition of local objects and local equivalences of \cite{hirsch}.
\begin{defn}\label{local}
Let $\Ch$ be a simplicial model category. Let $K$ be a set of objects in $\Ch$.
\begin{enumerate}
\item A map $f:A\to B$ is called a {\rm $K$-local equivalence} if for every element $X$ of $K$ the induced map of simplicial mapping spaces $f^{\ast}:\Map(B,X)\to \Map(A,X)$ is a weak equivalence of simplicial sets.
\item Let $\Kh$ denote the class of $K$-local equivalences. An object $X$ is called {\rm $K$-local} if it is $\Kh$-local,  i.e. if $X$ is fibrant and for every $K$-local equivalence $f:A\to B$ the induced
map $f^{\ast}:\Map(B,X) \to \Map(A,X)$ is a weak equivalence.
\end{enumerate}
\end{defn}
\begin{theorem}\label{leftloc}
Let $\Ch$ be a left proper fibrantly generated simplicial model category with all small limits. Let $K$ be a set of fibrant objects in $\Ch$ and let $\Kh$ be the class of $K$-local equivalences. 
\begin{enumerate}
\item The left Bousfield localization of $\Ch$ with respect to $\Kh$ exists, i.e.
there is a model category structure $\LKC$ on the underlying category $\Ch$ in which
\begin{enumerate}
\item the class of weak equivalences of $\LKC$ equals the class of $K$-local equivalences of $\Ch$,
\item the class of cofibrations of $\LKC$ equals the class of cofibrations of $\Ch$,
\item the class of fibrations of $\LKC$ is the class of maps with the right lifting property with respect
to those maps that are both cofibrations and $K$-local equivalences.
\end{enumerate}
\item The fibrant objects of $\LKC$ are the $K$-local objects of $\Ch$.
\item $\LKC$ is left proper. It is fibrantly generated if every object of $\Ch$
is cofibrant.
\item The simplicial structure of $\Ch $ gives $\LKC$ the structure of a simplicial model category.
\end{enumerate}
\end{theorem}
\begin{proof}
The proof is a dual version of the proof of Theorem 5.1.1 of \cite{hirsch}. But since we do not assume $\Ch$ to be cellular, we have to use a dual version of Lemma 2.5 of \cite{prospectra}, in order to avoid the use of Proposition 5.2.3 of \cite{hirsch}. One may look at \cite{pec} for more details.
\end{proof}

We return to spectra. Let $\Ch$ be a left proper fibrantly generated simplicial model category with all small limits and finite colimits. Let $T:\Ch \to \Ch$ be a left Quillen endofunctor on $\Ch$. Let $U:\Ch \to \Ch$ be its right adjoint.
\begin{defn}
A {\rm spectrum} $X$ is a sequence $(X_n)_{n\geq 0}$ of objects of $\Ch$ together
with structure maps $\sigma :TX_n \to X_{n+1}$ for all $n\geq 0$. A {\rm map of spectra}
$f:X\to Y$ is a collection of maps $f_n:X_n \to Y_n$ commuting with the structure maps.
We denote the category of spectra by $\SpC$.
\end{defn}
We begin by defining an intermediate strict structure, see \cite{bousfried} or \cite{hovey}. As usual we have the following functors:
\begin{defn}\label{FnRn}
Given $n\geq 0$, the {\rm evaluation functor} $\Ev_n:\SpC \to \Ch$ takes $X$ to $X_n$. It has a left adjoint $F_n:\Ch \to \SpC$ defined by $(F_nA)_m=T^{m-n}A$ if $m\geq n$ and $(F_nA)_m=\ast$ otherwise. The structure maps are the obvious ones.\\
The evaluation functor also has a right adjoint $R_n:\Ch \to \SpC$ defined by $(R_nA)_i=U^{n-i} A$ if $i\leq n$ and $(R_nA)_i=\ast$ otherwise. The structure map $TU^{n-i}A \to U^{n-i} A$ is adjoint to the identity map of $U^{n-i}$ when $i<n$.
\end{defn}
\begin{defn}\label{defstrict}
A map $f$ in $\SpC$ is a {\rm projective weak equivalence} (resp. {\rm projective fibration}) if each map $f_n$ is a weak equivalence (resp. fibration). A map $i$ is a {\rm projective cofibration} if it has the left lifting property with respect to all projective trivial fibrations. 
\end{defn}
The following proposition is proved in the standard way, \cite{hovey}, Proposition 1.14.
\begin{prop}\label{projcofs}
A map  $i:A \to B$ in $\SpC$ is a projective (trivial) cofibration if and only if $i_0:A_0 \to B_0$ and the induced maps  $j_n:A_n \amalg_{TA_{n-1}}TB_{n-1}$ for $n\geq 1$ are (trivial) cofibrations in $\Ch$. 
\end{prop}
We will show that the projective structure is in fact a fibrantly generated model structure. We denote the set of generating fibrations of $\Ch$ by $P$ and the set of generating trivial fibrations by $Q$. Inspired by Proposition \ref{projcofs}, we set 
$$\tP:=
\{g_n: R_nR \to R_nS\times_{R_{n-1}US}R_{n-1}UR, ~ \mbox{for all} ~ f:R\to S ~\mathrm{in} ~ P, n\geq 0\}$$ 
and 
$$\tQ:=
\{g_n: R_nR \to R_nS\times_{R_{n-1}US}R_{n-1}UR, ~ \mbox{for all} ~ f:R\to S ~\mathrm{in} ~ Q, n\geq 0\}$$ 
where $g_n$ is the map induced by the commutative diagram
\begin{equation}\label{diagRn}
\begin{array}{ccc}
R_nR& \to & R_nS \\
\downarrow &  & \downarrow \\
R_{n-1}UR & \to & R_{n-1}US.
\end{array}
\end{equation}
\begin{lemma}\label{levelfibs}
If $f:X \to Y$ is a fibration (resp. trivial fibration) in $\Ch$, then the maps 
$g_n: R_nX \to R_nY\times_{R_{n-1}UY}R_{n-1}UX$ are fibrations (resp. trivial fibrations) in $\SpC$.
\end{lemma}
\begin{proof}
Since $R_n$ is defined via the right Quillen functor $U$ which preserves fibrations and trivial fibrations, it is clear that the map $(R_nf)_i=U^{n-i}(f)$ on the $i$-th level is a fibration (resp. trivial fibration). In  diagram (\ref{diagRn}), the maps on the $i$-th level are either identities or equal to $U^{n-i}(f)$. Hence the induced map $g_n$ is also a fibration (resp. trivial fibration).
\end{proof}
\begin{theorem}\label{thmstrict}
Let $\Ch$ be a left proper fibrantly generated simplicial model category with all small limits.
The projective weak equivalences, projective fibrations and projective cofibrations define a left proper
fibrantly generated simplicial model structure on $\SpC$ with set of generating fibrations $\tP$ and set of generating trivial fibrations $\tQ$.
\end{theorem}
\begin{proof}
This intermediate result can be proven in essentially the same way as Theorem 1.14 in \cite{hovey}. The hard point is the factorization axiom which follows from a cosmall object argument and the fact that $\Ch$ is fibrantly generated. Since we use a cosmall instead of small object argument the existence of small limits suffices.  
\end{proof}

It remains to modify this structure in order to get a stable structure, i.e. one in which the prolongation of $T$ is a Quillen equivalence. We will do this by applying the localization Theorem \ref{leftloc} to the projective model structure on spectra. We want the stable weak equivalences to be the maps that induce isomorphisms on all generalized cohomology theories. A generalized cohomology theory is represented by the analogue of an $\Omega$-spectrum.
\begin{defn}
A spectrum $E\in \SpC$ is defined to be a {\rm $U$-spectrum} if each $E_n$  is fibrant and the adjoint structure maps $E_n \to UE_{n+1}$ are weak equivalences for all $n\geq 0$.
\end{defn}
By Corollary 9.7.5 in \cite{hirsch}, since each $E_n$ is fibrant and since the right Quillen functor $U$ preserves fibrations, we know that $E_n \to UE_{n+1}$ is a weak equivalence in $\Ch$ if and only if the induced map $\Map (A,E_n) \to \Map (A, UE_{n+1})$ of mapping spaces is a weak equivalence in $\Sh$ for every cofibrant object $A$. By adjunction this is equivalent to $\Map (F_nA,E) \to \Map (F_{n+1}TA, E)$ being a weak equivalence in $\Sh$ for every cofibrant object $A$. So we have to define  the set $K$ to consist of all $U$-spectra such that the maps $F_{n+1}TA \to F_nA$ adjoint to the identity map of $TA$ are the $K$-local equivalences. By the fact that the projective model structure on $\SpC$ is fibrantly generated, left proper and simplicial, we may deduce the following result from Theorem \ref{leftloc}.
\begin{theorem}\label{stablestructure}
Let $\Ch$ be a left proper fibrantly generated simplicial model category with all small limits. Let $K$ be the set of $U$-spectra. There is a {\rm stable model structure} on $\SpC$ which is defined to be the $K$-localized model structure $L_{\Kh}\SpC$ of Theorem \ref{leftloc} where $\Kh$ is the class of all $K$-local equivalences.
\end{theorem}
We apply this to the category $\Ch=\hShp$ and $T=S^1\wedge \cdot$. By Theorem \ref{fibgen} its model structure is simplicial, left proper and fibrantly generated. Note that $S^1\wedge \cdot$ is a left Quillen endofunctor since it takes monomorphisms to monomorphisms and preserves cohomological equivalences by the suspension axiom. Its right adjoint is $\hat{\Omega}$.
\begin{defn}\label{defspectrum}
A {\rm profinite spectrum} $X$ consists of a sequence $X_n \in \hShp$ of pointed profinite spaces for $n\geq0$ and maps $\sigma_n:S^1 \wedge X_n \to X_{n+1}$ in $\hShp$.\\ 
A {\rm morphism} $f:X \to Y$ of spectra consists of maps $f_n:X_n \to Y_n$ in $\hShp$ for $n\geq0$
such that $\sigma_n(1\wedge f_n)=f_{n+1}\sigma_n$.\\
We denote by $\hSp$ the corresponding category and call it the category of {\rm profinite spectra}.
\end{defn}
\begin{cor}\label{thmprostable}
There is a {\rm stable model structure on $\hSp$} for which the prolongation $S^1 \wedge \cdot:\hSp \to \hSp$ is a Quillen equivalence.\\
In particular, the {\rm stable equivalences} are the maps that induce an isomorphism on
all generalized cohomology theories, represented by profinite $\hat{\Omega}$-spectra;
the {\rm stable cofibrations} are the maps $i:A \to B$ such that $i_0$ and the induced maps $j_n: A_n \amalg_{S^1\wedge A_{n-1}}S^1 \wedge B_{n-1} \to B_n$ are monomorphisms for all $n$;
the {\rm stable fibrations} are the maps with the right lifting property
with respect to all maps that are both stable equivalences and stable cofibrations.
\end{cor}
\begin{remark}
There is also a stable model structure on symmetric profinite spectra in the sense of \cite{symmspectra}, see \cite{pec}.
\end{remark}
\paragraph{Profinite completion of spectra}
Let $\Sp$ be the stable model structure of simplicial spectra defined in \cite{bousfried}. Let $\compl :\Sp \to \hSp$ be the profinite completion  applied levelwise that takes the spectrum $X$ to the profinite spectrum $\hat{X}$ whose structure maps are given by
$$S^1 \wedge \hat{X}_n \cong \hat{S^1 \wedge X_n} \stackrel{\hat{\sigma}}{\to} \hat{X}_{n+1}.$$
Let $|\cdot|: \hSp \to \Sp$ be the levelwise applied forgetful functor. Since the two functors $\compl$ and $|\cdot|$ on $\hShp$ commute with smash products, they form an adjoint pair of functors. The following result may be easily deduced from \cite{ensprofin}, \S 2,Proposition 1.
\begin{prop}\label{quillenpair}
The functor $\compl :\Sp \to \hSp$ preserves stable equivalences and cofibrations.\\
The functor $|\cdot|: \hSp \to \Sp$ preserves fibrations and stable equivalences between
fibrant objects.\\
In particular, $\compl$ induces a functor on the homotopy categories and the pair $(\compl , |\cdot|)$ is a Quillen pair of adjoint functors.
\end{prop}
Similar versions of the following facts have also been proved by Dehon \cite{dehon}. We deduce the following result from Proposition \ref{piofZl-completion}. 
\begin{prop}\label{stablepiofcompletion}
Let $E$ be a $(-1)$-connected spectrum and suppose that the $\Zl$-cohomology of each $E_n$, $n \geq 1$, is finite dimensional in each degree. Then the stable homotopy groups of the profinite completion $\hat{E}$ are given by the following isomorphism for all $n$
$$\pi_n \hat{E} \cong \hat{\pi_n E}^{\ell} \cong \Z_{\ell} \otimes_{\Z} \pi_n E.$$
\end{prop}
\begin{example}\label{piofMU}
1. Let $MU$ be the simplicial spectrum representing complex cobordism and let $MU/\ell^{\nu}$ be the cofibre spectrum defined by multiplication with $\ell^{\nu}$. For the profinite completions there are isomorphisms
$$\pi_{2\ast}\hMU \cong  \Z_{\ell} \otimes_{\Z} \pi_{2\ast} MU ~\mbox{and}~\pi_{2\ast}\hMU/\ell^{\nu}\cong \Zln \otimes \pi_{2\ast}MU,$$
where $\pi_{2\ast}MU$ is a polynomial ring in infinitely many variables, cf. \cite{adams}.\\
2. Let $\hKU$ be the profinite completion of the $\Omega$-spectrum representing
complex $K$-theory: 
$\hKU_{2i}= \hat{BU\times  \Z}$ and $\hKU_{2i+1}= \hat{U}$
for all $i\geq 0$. Although $KU$ does not satisfy the hypothesis of the proposition, we get by Proposition \ref{piofZl-completion}
$\pi_{2i}(\hKU)\cong \Z_{\ell}$ and $\pi_{2i+1}(\hKU)= 0$ for all $i$. 
\end{example}
\paragraph{Generalized cohomology theories on profinite spaces}
It is the main feature of the stable homotopy category of profinite spectra that it provides a general setting for cohomology theories on profinite spaces. We define generalized cohomology theories on $\hSh$ to be the functors represented by profinite spectra. Profinitely completed cohomology theories have already been studied by Dehon in \cite{dehon}.\\
Let $E$ be a spectrum in $\hSp$, we set $E^n(X):= \Hom_{\hSHh}(X,E[n])$, where $\Hom_{\hSHh}(X,E[n])$ denotes the set of maps that lower the dimension by $n$, and call this the $n$-th cohomology group of $X$ with values in $E$. We set $E^{\ast}(X):=\bigoplus_n E^n(X)$.
For a pointed profinite space $X$ we define its $n$-th cohomology groups with values in $E$ by
\begin{equation}\label{defgencoho}
E^n(X):= \Hom_{\hSHh}(\Sigma ^\infty (X),E[n]).
\end{equation}
For a pair $(X,A)$ of profinite spaces we define the relative cohomology by 
$E^n(X,A):= E^n(X/A)$.
%
\paragraph{The Atiyah-Hirzebruch spectral sequence}
Let $X$ be a profinite space. For every integer $p$, we denote by $\sk_pX$ the profinite subspace of $X$ which is generated by the simplices of degree less or equal $p$. We call $\sk_pX$ the $p$-skeleton of $X$. For every $k$ the set of $k$-simplices of $\sk_pX$ is closed in $X_k$ such that $\sk_pX$ is a simplicial profinite subset of $X$, cf. \cite{dehon}. Let $E$ be a profinite spectrum. We consider the skeletal filtration $\sk_0X \subset \ldots \subset \sk_pX \subset \sk_{p+1}X \subset \ldots \subset X$. It yields a filtration on $E^{\ast}X$ defined by $F^pE^{\ast}X:=\mathrm{Ker}(E^{\ast}X \to E^{\ast}\sk_{p-1}X)$.\\
Since the coefficient groups $E^q$ have a natural profinite structure, we may consider continuous cohomology with coefficients in $E^q$. The proof of the following theorem is essentially the one given by Adams in \cite{adams}. Dehon has translated it to the profinite setting in \cite{dehon}, Proposition 2.1.9. Although Dehon proves the following assertion only for special cohomology theories, the proof works in the category of profinite spectra for any profinite cohomology theory, see \cite{pec}. 
\begin{prop}\label{ahss}
For any profinite spectrum $E$ and for any profinite space $X$ there is a spectral sequence $\{E_r^{p,q}\}$ with $E_2^{p,q}\cong H^p(X;E^q) \Longrightarrow E^{p+q}(X)$.
The spectral sequence converges strongly, in the sense of \cite{adams} III, \S 8.2, to the graded term $E_{\infty}^{p,q}=F^pE^{\ast}X/F^{p+1}E^{\ast}X$ of the filtration on $E^{\ast}X$ if ${\lim_r}^1E_r^{p,q}=0$. In particular, the spectral sequence converges if $H^p(X;E^q)=0$ for $p$ sufficiently large or if $H^p(X;E^q)$ is finite for all $p$.  
We call it the {\em profinite Atiyah-Hirzebruch spectral sequence}. 
\end{prop} 
\paragraph{Comparison with generalized cohomology theories of pro-spectra}
As a corollary, we compare the cohomology of profinite spectra with cohomology theories of pro-spectra. Isaksen constructs in \cite{gencohomofpro} a stable model structure on the category of pro-spectra. If $E$ is a pro-spectrum then the $r$-th cohomology $E_{\pro}^r(X)$ of a pro-spectrum $X$ with coefficients in $E$ is the set $[X,E]^{-r}_{\pro}$ of maps that lower the degree by $r$ in the stable homotopy category of pro-spectra. In addition, for two pro-spectra $X$ and $E$ he shows the existence of an Atiyah-Hirzebruch spectral sequence $E_2^{p,q}=H^{-p}(X;\pi_{-q}E)\Rightarrow [X,E]^{p+q}_{\pro}$ which is in particular convergent if $E$ is a constant pro-spectrum and $X$ is the suspension spectrum of a finite dimensional pro-space. Taking completion and homotopy limits defines a functor $\compl :  \pro-\Sp \to \hSp$. One can check that it also induces maps $E_{\pro}^{\ast}(X) \to \hat{E}^{\ast}(\hat{X})$, which yield morphisms of Atiyah-Hirzebruch spectral sequences. Since the $E_2$-terms agree, we get
\begin{cor}\label{proMUcomparison}
Let $\{X_s\}_{s \in I}$ be a finite dimensional pro-space. Then there is an isomorphism 
$$MU^{\ast}_{\pro}(\{X_s\};\Zln)\cong \hMU^{\ast}(\hat{X};\Zln).$$ 
\end{cor} 

\section{Stable \'etale realization on the $\A^1$-homotopy category}
The construction of the \'etale topological type functor $\Et$ from locally noetherian schemes to pro-spaces is due to Artin-Mazur and Friedlander. The construction of the $\A^1$-homotopy category of schemes gave rise to the question if this functor may be enlarged to the category of spaces. This has been answered independently by Isaksen and Schmidt. The latter one constructed in \cite{schmidt} a geometric functor to the category of pro-objects in $\Hh$. Isaksen gave a rigid, but less intuitive construction in \cite{a1real}. We follow Isaksen's approach. The first step in this direction was the construction of a model structure on $\pro-\Sh$. 
\paragraph{The functor $\hEt$}
Let $\pro-\Sh$ be the category whose objects are cofiltered diagrams $X(-):I \to \Sh$ and its morphisms are defined by 
$$\Hom_{\pro-\Sh}(X(-),Y(-)):=\lim_{t\in J}\colim_{s\in I} \Hom_{\Sh}(X(s),Y(t)).$$
The cohomology with $\Zl$-coefficients is defined to be
$$H^{\ast}(X(-),\Zl):=\colim_s H^{\ast}(X(s),\Zl).$$ 
Isaksen has constructed several model structures on $\pro-\Sh$; in particular, the $\Zl$-cohomological model structure of \cite{compofpro} in which the weak equivalences are morphisms inducing isomorphisms in $\Zl$-cohomology and the cofibrations are levelwise monomorphisms. \\
We define a completion functor $\compl:\pro-\Sh \to \hSh$ as the composite of two functors. First we apply $\compl :\Sh \to \hSh$ levelwise, then we take the limit in $\hSh$ of the underlying diagram. The next lemma shows that the functor $\compl$ respects the model structures on $\pro-\Sh$ and $\hSh$. Its  proof is clear using Lemme 1.1.1 of \cite{dehon}. Using the functor $\hSh \stackrel{G}{\longrightarrow} \pro-s\Fh \stackrel{i}{\longrightarrow} \pro-\Sh$ of Isaksen's it follows as in \cite{compofpro} from this lemma that the homotopy categories $\Ho(\pro-\Sh)$ and $\hHh$ are in fact equivalent via completion.
\begin{lemma}\label{Z/ell-comparison}
1. Let $X \in \pro-\Sh$ be a pro-simplicial set. Then we have a natural isomorphism of cohomology groups
$$H^{\ast}(X;\Zl) \stackrel{\cong}{\longrightarrow} H^{\ast}(\hat{X};\Zl).$$
2. A morphism $f:X\to Y$ of pro-simplicial sets induces an isomorphism in $\Zl$-cohomology if and 
only if the morphism $\hat{f}:\hat{X}\to \hat{Y}$ in $\hSh$  induces an isomorphism in continuous $\Zl$-cohomology. \\
3. The functor $\compl :\pro-\Sh \to \hSh$ preserves monomorphisms.
\end{lemma}
Now we turn to the applications in algebraic geometry. We refer the reader to \cite{fried} and \cite{a1real} for a detailed discussion of the category of rigid hypercoverings and rigid pullbacks. The following definition is taken from \cite{a1real}:
\begin{defn} 
Let $X$ be a locally noetherian scheme. The {\rm \'etale topological type of $X$} is defined to be the pro-simplicial set
$$\Et X := \mathrm{Re} \circ \pi :HRR(X) \to \Sh$$
sending a rigid hypercovering $U_{\cdot}$ of $X$ to the simplicial set of connected components of $U_{\cdot}$. If $f: X \to Y$ is a map of locally noetherian schemes, then the strict map $\Et f:\Et X \to \Et Y$ is given by the functor $f^{\ast}:HRR(Y) \to HRR(X)$ and the natural transformation $\Et X \circ \Et f \to \Et Y$.
\end{defn}
Isaksen uses the insight of Dugger \cite{dugger} that one can construct the unstable $\A^1$-homotopy category in a universal way. Starting from an almost arbitrary category $\Ch$, Dugger constructs an enlargement of $\Ch$ that carries a model structure and is universal for this property. Isaksen extends in \cite{a1real} the functor $\Et$ via this general method to the $\A^1$-homotopy category. The idea is that $\Et X$ should be the above $\Et X$ on a representable presheaf $X$ and should preserve colimits and the simplicial structure. For our purpose, we would like to define $\hEt$ directly in the way Dugger suggests. But the problem is that $\compl:\pro-\Sh \to \hSh$ is not a left adjoint functor and does not preserve all small colimits. Therefore, we define $\hEt$ to be the composition of $\Et$ followed by completion.\\ 
Let $\Sm/k$ be the category of smooth quasi-projective schemes of finite type over $k$ and $\Delta^{\mathrm{op}}\mathrm{PreShv}(\Sm/k)$ the category of simplicial presheaves on $\Sm/k$. We make the following definition whose first part is due to Isaksen \cite{a1real}:
\begin{defn}\label{defet}
If $X$ is a representable presheaf, then $\Et X$ is the \'etale topological type of $X$.
If $P$ is a discrete presheaf, i.e. $P$ is
just a presheaf of sets, then $P$ can be written as a colimit $\colim_i X_i$ of representables and we define $\Et P:=\colim_i \Et X_i$. Finally, an arbitrary simplicial presheaf can be written as the coequalizer
of the diagram
$$\coprod_{[m]\to[n]} P_m \otimes \Delta^n \rightrightarrows \coprod_{[n]} P_n \otimes \Delta^n,$$
where each $P_n$ is discrete. Define $\Et P$ to be the coequalizer of the diagram
$$\coprod_{[m]\to[n]} \Et P_m \otimes \Delta^n \rightrightarrows \coprod_{[n]} \Et P_n \otimes \Delta^n.$$
We define the {\rm profinite \'etale topological type functor} $\hEt$ to be the composition of $\Et$ and the profinite completion functor $\pro-\Sh \to \hSh$:
$$\hEt := \compl \circ \Et : \Delta^{\mathrm{op}}\mathrm{PreShv}(\Sm/k) \to \hSh.$$
\end{defn}
For computations the following remark is crucial. 
\begin{remark}\label{remarketalecohom}
Let $M$ be a locally constant sheaf on the locally noetherian scheme $X$. By \cite{fried}, Proposition 5.9, we know that there is an isomorphism for \'etale cohomology $H^{\ast}_{\et}(X;M) \cong H^{\ast}(\Et X,M)$ where $H^{\ast}(\Et X;M)$ denotes the cohomology of a pro-simplicial set $\Et X$ with coefficients in the local coefficient system $M$ corresponding to the sheaf $M$. For a finite abelian group $\pi$, we have in addition a natural isomorphism by Lemma \ref{Z/ell-comparison}: $H^{\ast}(Z;\pi) \cong H^{\ast}(\hat{Z};\pi)$ for every pro-simplicial set $Z$. Hence we get as well $$H^{\ast}_{\et}(X;\pi) \cong H^{\ast}(\hEt X,\pi)$$ for every locally noetherian scheme $X$ and every finite abelian group $\pi$.\\
Furthermore, the pro-fundamental group $\pi_1(\Et X)$ is isomorphic as a profinite group to the \'etale fundamental group $\pi_1^{\et}(X)$ of $X$. Hence $\pi_1(\hEt X)$ is equal to $\pi_1^{\et}(X)^{\wedge \ell}$ by Proposition \ref{piofZl-completion}.
\end{remark}
\paragraph{Examples}
1. Let $R$ be a strict local henselian ring, i.e. a local henselian ring with separably closed residue field.
Then $\Spec R$ has no nontrivial \'etale covers and the \'etale topological type of $\Spec R$ is a contractible space.\\
2. Let $k$ be a separably closed field with $\mathrm{char}(k)\neq \ell$. The space $\Ge_m$ is connected and its $\ell$-completed \'etale fundamental group is $\Z_{\ell}$, the $\ell$-adic integers. This implies $\hEt \Ge_m \cong K(\Z_{\ell},1)$ in $\hSh$, where $K(\Z_{\ell},1)$ is naturally a profinite space since $\Z_{\ell}$ is a profinite group.\\
3. Let $k$ be a separably closed field with $\mathrm{char}(k)\neq \ell$. Let $\Pro_k^1$ be the projective line over $k$. Since $\Pro_k^1$ is connected and its \'etale fundamental group $\pi^{\et}_1 (\Pro_k^1)=0$ is trivial, see e.g.~\cite{milne} I, Example 5.2 f), $\hEt \Pro_k^1$ is a simply connected space. Apart from $H^0$ its only nonzero \'etale cohomology group is $H^2_{\et}(\Pro^1_k ,\Z/\ell) \cong \Zl$, with a chosen isomorphism $\Zl \to \Zl(1)$. Hence $\hEt \Pro_k^1$ is isomorphic in $\hSh$ to the simplicial finite set $S^2$.\\
4. Let $k=\F_q$ be a finite field with $\mathrm{char}(k)=p\neq \ell$. The \'etale topological type of $k$ is isomorphic to $S^1$ in $\hHh$. For, it is connected and its $\ell$-completed fundamental group is the $\ell$-completion of the absolute Galois group of $k$, i.e. $\pi^{\ell}_{1,\et}(k)=\Z_{\ell}$. Since the cohomology groups $H^i(k;\Zl)$ vanish for $i>1$, $\hEt k$ is a space of dimension one and weakly equivalent to $S^1$. 
\subsection{Profinite \'etale realization of motivic spaces}
We consider as in \cite{dugger} the category $U(\Sm/k):=\Delta^{\mathrm{op}}\mathrm{PreShv}(\Sm/k)$ of simplicial presheaves on $\Sm/k$ with the projective model structure: the weak equivalences (fibrations) are objectwise weak equivalences (fibrations) of simplicial sets, the cofibrations are the maps having the left lifting property with respect to all trivial fibrations.
Then one takes the left Bousfield localization of this model structure at the set $S$ of maps:
\begin{enumerate}
\item\label{rel1} for every finite collection $\{X^a\}$ of schemes with disjoint union $X$, the
map $\coprod X^a \to X$ from the coproduct of the presheaves represented by $X^a$
to the presheaf represented by $X$;
\item\label{rel2} every Nisnevich (\'etale) hypercover $U_{\cdot} \to X$;
\item\label{rel3} $X\times \A^1 \to X$ for every scheme $X$.
\end{enumerate}
We call this the {\em Nisnevich (\'etale) $\A^1$-local projective model structure}
according to \cite{a1real}, and denote it by $\LU$ (resp. $\LeU$).\\
Proposition 8.1 of \cite{dugger} states that $\LU$ is Quillen equivalent to the Nisnevich $\A^1$-localized model category $\MVh$ of \cite{mv} and the analogue holds for the \'etale case. For an \'etale realization over a field of positive characteristic, we have to consider the following point. Since the \'etale fundamental group of the affine line $\A^1_k$ over a field $k$ of positive characteristic is non-trivial, one has to complete away from the characteristic of $k$ in this case, using the fact that the projection $X\times \A^1 \to X$ induces an isomorphism in \'etale cohomology $H_{\et}^{\ast}(X;\Zl) \to H_{\et}^{\ast}(X\times \A^1;\Zl)$ for every prime $\ell \neq p$. This means that we have to consider the $\Zl$-cohomological model structure on pro-$\Sh$, resp. $\hSh$. The following theorem is due to Isaksen \cite{a1real}:
\begin{theorem}\label{LhEt}
Let $\ell$ be a prime different from the characteristic of $k$. With respect to the Nisnevich (\'etale) $\A^1$-local projective model structure on simplicial presheaves on $\Sm/k$, the functor $\hEt$ induces a functor $\LhEt$ from the Nisnevich (\'etale) $\A^1$-homotopy category to the $\Zl$-cohomological homotopy category of $\hSh$.
In particular, $\LhEt X$ is just $\hEt X$ for every scheme in $\Sm/k$, and hence $\hEt$ preserves $\A^1$-weak equivalences between smooth schemes over $k$.
\end{theorem}
\begin{proof}
Since every Nisnevich hypercover is also an \'etale hypercover it suffices to check this for the \'etale case. The first part of the proof is the one of \cite{a1real}, Theorem 6.\\
By Lemma \ref{Z/ell-comparison} completion preserves weak equivalences and cofibrations. Hence the composition $\hEt$ sends weak equivalences between cofibrant objects into weak equivalences and the total left derived functor $\LhEt: \Hh _{\A^1}^{\et}(k) \to \hHh$ from the \'etale $\A^1$-homotopy category $\Hh _{\A^1}^{\et}(k)$ of schemes over $k$ of \cite{mv} exists.
\end{proof}
\subsection{Etale realization of motivic spectra}
We extend the results of the previous section to the stable $\A^1$-homotopy category. This \'etale realization of the stable motivic category is the technical key point for the construction of a transformation from algebraic cobordism given by the $MGL$-spectrum to the profinite \'etale cobordism of the next section. One of the reasons why we consider the model $\hSh$ of $\Ho(\hSh)$, instead of $\pro-\Sh$, is that it seems to be easier to extend the functor $\hEt$ to the category of profinite spectra rather than pro-spectra. I am especially grateful to Fabien Morel for helpful discussions on the problems in this section.\\ 
Let $k$ be the base field of characteristic different from $\ell$ and let $\ok$ be its separable closure. For the category $\mathrm{Sp}^{\Pro^1}(k)$ of motivic $\Pro^1_k$-spectra over $k$, in the sense of \cite{pi0} for example, we have to consider presheaves $\Xh$ pointed by a morphism $\Spec k \to \Xh$. This forces us to consider also the category $\hShpk$ of pointed profinite spaces over $\hEt k$. Its objects $(X,p,s)$ are pointed profinite spaces $X$ together with a projection morphism  $p:X \to \hEt k$ and a section morphism $s:\hEt k \to X$. The morphisms in this category are commutative diagrams in the obvious sense. The \'etale realization of a pointed presheaf $\Xh$ is naturally an object of $\hShpk$ via the images of the projection and section morphisms of $\Xh$. For the contractible space $\hEt \ok$, we choose and fix a base point. Since all \'etale covers of $\Spec \ok$ are trivial, any choice is ok. The point is that this allows us to consider every $\hEt \Xh$ as a pointed space compatible with base extensions of $k$ in the following way. The space $\hEt \Xh$ is pointed by the composite $\ast \to \hEt \ok \to \hEt k \to \Xh$. Furthermore, via the canonical maps $X \to \ast \to \hEt \ok \to \hEt k$ and $\hEt k \to \ast \to X$ we may view every space $X$ in $\hShp$ as an object in $\hShpk$. In particular, the $2$-sphere $S^2$ is naturally an object in $\hShpk$. We consider the usual model structure on $\hShpk$ where weak equivalences (resp. cofibrations, fibrations) are those maps which are $\Zl$-weak equivalences (resp. cofibrations, fibrations) in $\hSh$ after forgetting the projection and section maps. This model structure on $\hShpk$ is left proper and fibrantly generated and we construct a stable model structure on the category $\mathrm{Sp}(\hShpk, S^2 \wedge \cdot)$ of profinite $S^2$-spectra over $\hEt k$ exactly in the same way as for $\mathrm{Sp}(\hShp, S^2 \wedge \cdot)$. Its homotopy category will be denoted by $\hSHhtk$. 
\paragraph{An intermediate category} The problem for the stable motivic version is that $\Et$ and hence also $\hEt$ do not commute with products in general. We have to construct an intermediate category and show by a zig-zag of functors that we get a functor on the homotopy level. However, the projections to each factor induce a canonical map $$\hEt(\Pro^1_k \wedge_k X) \to \hEt(\Pro^1_k) \wedge_{\hEt k} \hEt(X).$$

\begin{lemma}\label{Kunneth}
For every pointed presheaf $\Xh$ on $\Sm/k$, the sequence of canonical maps in $\hShpk$
$$S^2 \wedge \hEt \Xh \stackrel{\simeq}{\longrightarrow} \hEt (\Pro^1_k) \wedge_{\hEt k} \hEt(\Xh)
\stackrel{\simeq}{\longleftarrow} \hEt(\Pro^1_k \wedge_k \Xh)$$
is a sequence of $\Zl$-weak equivalences.
\end{lemma}
\begin{proof}
For $X \in \Sm/k$, the assertion may easily be deduced from the projective bundle formula for \'etale cohomology and the isomorphism $\hEt \Pro^1_{\ok} \cong S^2$ in $\hShp$.\\
If $\Xh$ denotes a presheaf on $\Sm/k$, $\Xh$ is isomorphic to the colimit of representable presheaves $\Xh =\colim_s X_s$. Since each $X_s$ is a smooth scheme over $k$, the \'etale cohomology groups $H^i(\Et X_s;\Zl) = H_{\et}^i(X_s;\Zl)$ are finite $\Zl$-vector spaces in each degree. Hence, since $\Et $ commutes with colimits, we get $H^{\ast}(\Et \Xh;\Zl) \cong \lim_s H^{\ast}(\Et X_s;\Zl)$, as the $\lim^1$ of these finite groups vanishes, and the limit over all $s$ of $H^{\ast}(\Et X_s;\Zl)$ commutes with the functor $H^{\ast}(\Et \Pro^1_k;\Zl) \otimes_{H^{\ast}(\Et k;\Zl)}-$. 
\end{proof}

Hence if $\sigma_n:\Pro^1_k \wedge_k E_n \to E_{n+1}$ is the structure map of a motivic $\Pro^1_k$-spectrum, then $\hEt$ yields a sequence of maps 
\begin{equation}\label{structuremap}
S^2 \wedge \hEt E_n \stackrel{\simeq}{\longrightarrow} \hEt (\Pro^1_k) \wedge_{\hEt k} \hEt( E_n)
\stackrel{\simeq}{\longleftarrow} \hEt(\Pro^1_k \wedge_k E_n) \stackrel{\hEt \sigma_n}{\longrightarrow}
E_{n+1}
\end{equation}
where the first two maps are weak equivalences in $\hShpk$. Since there is no natural inverse map $\hEt (\Pro^1_k) \wedge_{\hEt k} \hEt(\Xh)
\to \hEt(\Pro^1_k \wedge_k \Xh)$ in $\hShpk$, we may only construct an \'etale realization on the level of homotopy categories. Therefore, we consider an intermediate category $\Ch$ and deduce from a zig-zag of functors 
$$\mathrm{Sp}^{\Pro^1}(k) \stackrel{\tilde{\hEt}}{\to} \Chk \stackrel{i}{\hookleftarrow} \mathrm{Sp}(\hShk, S^2 \wedge \cdot)$$
the existence of a stable realization functor $\SHh(k) \to \hSHhtk$.\\
In view of Lemma \ref{Kunneth}, it is natural to consider the following definition. The objects of the category $\Chk$ are sequences
$$\{ F_n, F'_n, F''_n; S^2 \wedge F_n \stackrel{\simeq p_n}{\longrightarrow} F'_n \stackrel{\simeq q_n}{\longleftarrow} F''_n \stackrel{r_n}{\longrightarrow} F_{n+1} \}_{n \in \N}$$
where $F_n$, $F'_n$, $F''_n$ are pointed profinite spaces over $\hEt k$ and $p_n$, $q_n$ and $r_n$ are maps in $\hShpk$; furthermore the maps $p_n$ and $q_n$ are weak equivalences in $\hShp$.\\ 
The morphisms of $\Chk$ are levelwise morphisms of $\hShpk$ which make the obvious diagrams commutative, where the map $S^2 \wedge E_n \to S^2 \wedge F_n$ is the map induced by $E_n \to F_n$. The functor $\mathrm{Sp}(\hShpk, S^2 \wedge \cdot)\stackrel{i}{\hookrightarrow} \Chk$ denotes the full embedding which sends $\{ F_n, S^2 \wedge F_n \stackrel{\sigma_n}{\to} F_{n+1} \}$ to  $\{F_n, S^2 \wedge F_n, S^2 \wedge F_n; S^2 \wedge F_n \stackrel{\mathrm{id}}{\to} S^2 \wedge
F_n \stackrel{\mathrm{id}}{\leftarrow} S^2 \wedge F_n \stackrel{\sigma_n}{\to} F_{n+1} \}$.\\ 
On the other hand, we get a functor $\mathrm{Sp}^{\Pro^1}(k) \stackrel{\tilde{\hEt}}{\to} \Chk$, when we apply $\hEt$ levelwise.\\
We define a class $W$ of maps in $\Chk$ as the image of the stable equivalences of $\mathrm{Sp}(\hShpk, S^2 \wedge \cdot)$ under the embedding $i$. Since the maps in $W$ are the images of weak equivalences in a model structure and since $i$ is a full embedding, it follows that $W$ admits a calculus of fractions and we may form the localized category $\Ho(\Chk):=\Ch[W^{-1}]$. We call the maps in $W$ weak equivalences or stable equivalences, by abuse of notation. We will call a map in $W$ a level equivalence if it is in the image of the level equivalences of $\mathrm{Sp}(\hShpk, S^2 \wedge \cdot)$ under $i$. 
\begin{prop}\label{equivofcat}
The induced embedding $i: \hSHhtk \to \Ho(\Chk)$ is an equivalence of categories. 
\end{prop}
The corresponding inverse equivalence is denoted by $j: \Ho(\Chk) \to \hSHhtk$. Because of the necessary choices in the proof, one should note that $j$ is not quite constructible in practice.\\
\begin{proof}
Since $i$ is a full embedding, it suffices to show that for every $F \in \Ho(\Chk)$ there is a spectrum $E \in \hSHhtk$ such that $i(E)\cong F$ in $\Ho(\Chk)$. The crucial point is to construct the structure map of a spectrum from the given data of $F$.\\ 
Let $R$ be a fixed fibrant replacement functor in $\hShpk$. We consider the category $\mathrm{Sp}(\hShpk, RS^2 \wedge \cdot)$ as a Quillen equivalent model for $\hSHhtk$. Since $R$ commutes with products, applying $R$ on each level yields the following sequence
$$RS^2 \wedge RF_n \stackrel{\simeq Rp_n}{\longrightarrow} RF'_n \stackrel{\simeq Rq_n}{\longleftarrow} RF''_n \stackrel{Rr_n}{\longrightarrow} RF_{n+1}.$$
In addition, the functor $R$ can be chosen such that the map $Rq_n$ is a trivial fibration between fibrant and cofibrant objects, see Proposition 8.1.23 of \cite{hirsch}. The sequence we get is still isomorphic in $\Ho(\Chk)$ to the initial one since they are even level equivalent. By Proposition 9.6.4 of \cite{hirsch}, there is a right inverse $s_n$ of $Rq_n$ in $\hShpk$ such that $Rq_ns_n=\id_{RF'_n}$ and a homotopy $s_nRq_n \sim \id_{RF''_n}$. We denote by $E$ the resulting spectrum with structure maps $\sigma_n:=Rr_n\circ s_n\circ Rp_n$.\\ 
Now it is easy to check that $i(E)$ is isomorphic to $F$ in $\Ho(\Chk)$.
\end{proof}

\paragraph{The main result} We want to show that  $\thEtsp:\mathrm{Sp}^{\Pro^1}(k) \to \Chk$ has a total left derived functor. Therefore, we have to choose a good model for the stable motivic category $\SHh(k)$. By \cite{dugger} we know that $\LU$ is a left proper cellular simplicial model category. Hence, in order to be able to apply the methods of \cite{hovey}, we only need an appropriate left Quillen endofunctor. Unfortunately, as the referee pointed out, the map $\ast \to \Pro^1$ is not a projective cofibration in $\LU$. One should note that this problem does not appear when we use the injective model structure on presheaves, in which every monomorphism is a cofibration, see Jardine's work on motivic spectra \cite{jardine}. But $\Et$ is not a Quillen functor on the injective structure since there are too many injective cofibrations, see \cite{a1real} Remark 3.\\ 
Nevertheless, there are at least two possible approaches to solve this problem and we are indebted to the referee for calling our attention to them. The probably more elegant solution is to show that $\Et$ is a left Quillen functor on the flasque model structure on $\Delta^{\mathrm{op}}\mathrm{PreShv}(\Sm/k)$ constructed by Isaksen in \cite{flasque}, for which $\ast \to \Pro^1$ is a flasque cofibration. The easier and more obvious way is to factor $\ast \to \Pro^1$ as $\ast \to Q \to \Pro^1$, where $\ast \to Q$ is a projective cofibration and $p:Q \to \Pro^1$ is a trivial fibration in $\LU$. We follow the latter approach, which is in the spirit of \cite{motivicfunctors}, section 2.2, where the sphere $\A^1/\A^1-\{0\}$ is replaced by a cofibrant motivic space. Up to weak equivalence the choice of the motivic sphere is irrelevant for the stable motivic homotopy category, see \cite{hovey}, Theorem 5.7.  Since $\ast \to Q$ is a cofibration, the functor $Q \wedge \cdot:\LU \to \LU$ is a left Quillen functor and we are able to apply the methods of \cite{hovey}, sections 3 and 5. Furthermore, via the map $p:Q \to \Pro^1$, every $\Pro^1$-spectrum induces a $Q$-spectrum. 
\begin{prop}\label{p1model}
The canonical functors 
$$\mathrm{Sp}^{\Pro^1}(k) \to \SpMVP \to \SpLUQ$$
are Quillen equivalences.
\end{prop}
\begin{proof}
The first equivalence follows as in the proof of \cite{hovey}, Corollary 3.5, taking into account that cofibrations are sent to cofibrations and that the fibrant objects agree in both model structures.\\
The second equivalence follows from \cite{hovey}, Theorem 5.7, taking into account that $\LU$ is Quillen equivalent to $\MVh$ and that all objects in $\MVh$ are cofibrant. The functor $\MVh \to \LU$ is just the full embedding of Nisnevich sheaves into presheaves, its inverse is sheafification, see \cite{dugger}.
\end{proof}

Since $p$ is a trivial fibration between cofibrant objects, it is even a simplicial homotopy equivalence. Hence $\hEt p:\hEt Q \to \hEt \Pro^1$ is a weak equivalence in $\hShp$. The base change $\cdot \wedge_{\Spec k}\Spec \ok$ preserves homotopy equivalences. Hence, by the example above, $\hEt Q_{\ok}$ is a simply connected space whose cohomology is equal to the one of $S^2$. We conclude that $\hEt Q_{\ok}$ is even isomorphic to the simplicial finite set $S^2$  and to $\hEt \Pro^1_{\ok}$ in $\hShp$ and the sequence of maps in Lemma \ref{Kunneth} yields a sequence of $\Zl$-weak equivalences in $\hShp$ for every simplicial presheaf $\Xh$:
\begin{equation}\label{QKunneth}
S^2 \wedge \hEt \Xh \stackrel{\simeq}{\longrightarrow} \hEt (Q) \wedge_{\hEt k} \hEt(\Xh)
\stackrel{\simeq}{\longleftarrow} \hEt(Q \wedge_k \Xh).
\end{equation} 
Finally, we can prove that there is a stable \'etale realization functor on the $\A^1$-homotopy theory of schemes over an arbitrary base field of characteristic different from $\ell$.
\begin{theorem}\label{P1real}
The functor $\hEt$ induces an \'etale realization of the stable motivic homotopy category of $\Pro^1$-spectra: $$\LhEt: \SHh(k) \to \hSHhtk.$$
\end{theorem}
\begin{proof}
We use $ \SpLUQ$ as a model for $\SHh(k)$. The desired functor is defined to be the composite $\LhEt : \SHh(k) \stackrel{\LthEt}{\longrightarrow} \Ho(\Chk) \stackrel{j}{\longrightarrow} \hSHhtk$. For the existence of $\LthEt$ it remains to show that stable equivalences in $\SpLUQ$ are sent to isomorphisms in $\Ho(\Chk)$.\\
We know that $\thEtsp$ sends level equivalences between cofibrant objects in $\LU$ to weak equivalences in $\Chk$, since $\hEt$ sends weak equivalences between cofibrant objects to weak equivalences in $\hShpk$. Hence it induces a total left derived functor on the projective model structure of $\SpLUQ$.\\
We use the notation $\Sigma^Q_n:\LU \to \SpLUQ$ for the left adjoint to the $n$-th evaluation functor. It is given by $(\Sigma^Q_n \Xh)_m = Q^{m-n}\Xh$ if $m\geq n$ and $(\Sigma^Q_n \Xh)_m = \Spec k$ otherwise, where $Q^{m-n}\Xh$ denotes the smash product of $\Xh$ with $m-n$ copies of $Q$. We denote by $F_n:\hShpk \to \Chk$ the composition of the corresponding functor $\hat{\Sigma}_n:\hShpk \to \mathrm{Sp}(\hShpk,S^2\wedge \cdot)$ followed by the embedding $i$.\\ 
In order to show that there exists a derived functor on the stable structure it suffices to show that  $\hEtsp(\zeta^{\Xh}_n)$ is a stable equivalence for maps $\zeta^{\Xh}_n:\Sigma_{n+1}Q\wedge \Xh \to \Sigma_n \Xh$ in $ \SpLUQ$ for all cofibrant presheaves $\Xh \in \LU$, since these are the maps at which we localize to get the stable model structure, cf. \cite{hirsch} and \cite{hovey}.\\ 
We consider the commutative diagram in $\Ho(\Chk)$
$$\begin{array}{cccc} 
\zeta^{\hEt \Xh}_n: & F_{n+1}(S^2\wedge \hEt \Xh) & \stackrel{\cong}{\longrightarrow} & F_n(\hEt \Xh)\\
  & \updownarrow \cong &  & \parallel \\
       & F_{n+1}(\hEt Q \wedge_{\hEt k} \hEt \Xh))  & 
        & F_n(\hEt \Xh)\\
  & \updownarrow \cong &  & \parallel \\
       & F_{n+1}(\hEt (Q \wedge \Xh))  & 
        & F_n(\hEt \Xh)\\
 & \cong \updownarrow &  & \updownarrow \cong\\
\hEtsp\zeta^{\Xh}_n: & \hEtsp (\Sigma^Q_{n+1}(Q \wedge \Xh)) &
\longrightarrow & \hEtsp (\Sigma^Q_n \Xh).  
\end{array}$$
The upper and middle vertical isomorphisms on the left hand side are the obvious level equivalences deduced from the canonical sequence of weak equivalences (\ref{QKunneth}). The lower vertical isomorphisms are given by the following lemma applied to $\Xh$ and $Q\wedge \Xh$ respectively.
\begin{lemma}
For every integer $n\geq 0$ and every simplicial presheaf $\Xh$ there is an isomorphism 
$F_n(\hEt \Xh) \cong \hEtsp(\Sigma^Q_n \Xh)$ in $\Ho(\Chk)$.
\end{lemma}
\begin{proof}
We define an intermediate object $E^n \in \Chk$ given in degree $m\geq n$ by 
$$E^n_m= (\hEt Q)^{m-n}\wedge_{\hEt k} \hEt \Xh  , (E^n)_m'= (\hEt Q)^{m+1-n}\wedge_{\hEt k} \hEt \Xh, (E^n)_m''= (E^n)_m'$$ with the obvious structure maps induced by $S^2 \to \hEt Q$ respectively the identity; in degree $m<n$ it is defined by $E_m^n=(E^n)_m'=(E^n)_m''=\hEt k$ with identity maps and the map to the terminal object $\hEt k$ of $\hShpk$.\\
The object $E^n$ is defined such that there are canonical maps
$$F_n(\hEt \Xh) \stackrel{\alpha}{\longrightarrow} E^n \stackrel{\beta}{\longleftarrow} \hEtsp(\Sigma^Q_n \Xh)$$ induced in degree $m\geq n$ by the canonical weak equivalences of (\ref{QKunneth})
$$(S^2)^{m-n}\wedge \hEt \Xh \stackrel{\simeq}{\longrightarrow} (\hEt Q)^{m-n}\wedge_{\hEt k} \hEt \Xh 
\stackrel{\simeq}{\longleftarrow} \hEt (Q^{m-n}\wedge_k \Xh).$$
It follows that the maps $\alpha$ and $\beta$ are both level equivalences in $\Chk$. Hence $\alpha$ and $\beta$ are isomorphisms in $\Ho(\Chk)$. Their composition is the isomorphism $F_n(\hEt \Xh) \cong \hEtsp(\Sigma^Q_n \Xh)$.
\end{proof}
  
We deduce from the diagram that $\zeta^{\hEt \Xh}_n$ and $\hEtsp\zeta^{\Xh}_n$ differ only by an isomorphism in $\Ho(\Chk)$ given by level equivalences. Since $\zeta^{\hEt \Xh}_n$ is a stable equivalence in $\Chk$ by definition, the map $\hEtsp\zeta^{\Xh}_n$ is in fact an isomorphism in $\Ho(\Chk)$, which finishes the proof of the theorem.
\end{proof}

Let $MGL$ denote the motivic spectrum defined in \cite{a1hom} representing algebraic cobordism.
\begin{prop}\label{EtofMGL}
Let $k$ be a separably closed field. There is an isomorphism in $\hSHh_2$ 
$$\LhEt (MGL) \cong \hMU.$$ 
\end{prop}
\begin{proof}
Let $G_k(n,N)$ be the Grassmannian over $k$ and $G_k(n)$ be the colimit over $N$. By the Thom isomorphisms in \'etale and singular cohomology, it suffices to show that there is a weak equivalence between $\hEt (G_k(n))$ and $\hat{G}_{\C}(n)$, the profinite completion of the simplicial set corresponding to the complex Grassmannian manifold.\\
If $\car k >0$, let $R$ be the ring of Witt vectors of $k$, otherwise set $R=k$. By choosing a common embedding of $R$ and $\C$ into an algebraically closed field, Friedlander proved in \cite{etaleK}, 3.2.2, that there is a natural sequence of $\Zl$-weak equivalences in $\pro-\Sh_{\ast}$ 
$$\mathrm{Sing}(G_{\C}(n,N))  \to \Et(G_{\C}(n,N)) \to \Et(G_R(n,N)) \leftarrow \Et(G_k(n,N)).$$
We remark that weak equivalences in $\pro-\Hh_{\ast}$ in the sense of \cite{etaleK} correspond to $\Zl$-weak equivalences in $\pro-\Sh_{\ast}$ in the sense of \cite{compofpro}. By taking colimits with respect to $N$ we get a sequence of $\Zl$-weak equivalences in $\pro-\Sh_{\ast}$ 
$$\mathrm{Sing}(G_{\C}(n)) \to \Et(G_{\C}(n))  \to \Et(G_R(n)) \leftarrow \Et(G_k(n)).$$ 
Since $\Zl$-weak equivalences are preserved under completion, this shows that there is an isomorphism $\hat{G}_{\C}(n) \cong \hEt (G_k(n))$ in $\hHhp$.  
\end{proof}
\paragraph{Etale realization of $S^1$-spectra}
The \'etale realization of $S^1$-spectra is essentially simpler. The structure maps $\sigma_n:S^1_k \wedge_k E_n \to E_{n+1}$ induce by base extension canonical maps $S^1_{\ok} \wedge E_n \to S^1_k \wedge_k E_n \stackrel{\sigma_n}{\to} E_{n+1}$ in $\LU$. Together with the isomorphism $\hEt S^1_{\ok} = S^1$ in $\hShp$, we conclude that $\hEt$ defines a functor on $S^1$-spectra $\hEt : \mathrm{Sp}^{S^1}(k) \to \hSpk$. 
\begin{theorem}
The functor $\hEtsp : \mathrm{Sp}^{S^1}(k) \to \hSpk$ induces a functor $\LhEtsp : \SHhS \to \hSHhk$ on the stable $\A^1$-homotopy category of $S^1$-spectra. 
\end{theorem}

\section{Profinite \'etale cohomology theories}
\subsection{General theory}
Let $k$ be a fixed base field and let $\ell$ be a fixed prime different from the characteristic of $k$. The first main application of the previous discussion is a canonical setting for \'etale cohomology theories for schemes, in particular \'etale topological cobordism.
\begin{defn}\label{profinetalecohom}
Let $E\in \hSp$ be a profinite spectrum and let $X \in \Sm/k$. \\
1. We define the {\em profinite \'etale cohomology of $X$ in $E$} to be the profinite cohomology theory represented by $E$ applied to the profinite space $\hEt X$, i.e. $$E^n_{\et}(X):=E^n(\hEt X)=\Hom_{\hSHh}(\Sigma^{\infty}(\hEt X), E[n]),$$
where we add a base-point if $X$ is not already pointed.\\ 
2. We define relative \'etale cohomology groups $E^n_{\et}(X,U)$ for an open subscheme $U\subset X$ by $E^n_{\et}(X,U):= E^n(\hEt (X)/\hEt (U))$.\\
3. Similarly, if $\Xh$ is a simplicial presheaf on $\Sm/k$, we define the \'etale cohomology group of $\Xh$ to be $E_{\et}^n(\Xh):=\Hom_{\hSHh}(\Sigma^{\infty}(\hEt \Xh), E[n])$.
\end{defn}
The following fact may be shown in the standard way using Theorem 2.10 of \cite{a1real}.
\begin{prop}\label{mayervietoris}
Let $U \stackrel{i}{\hookrightarrow} X \stackrel{p}{\leftarrow} V$ induce an elementary distinguished square, \cite{a1hom}. Let $E$ be a profinite spectrum. Then there is a Mayer-Vietoris long exact sequence of graded groups
$$\cdots \to E^n_{\et}(X)\to E^n_{\et}(U)\oplus E^n_{\et}(V)\to E^n_{\et}(U\times_X V) \to E^{n+1}_{\et}(X)\to \cdots.$$  
\end{prop}
\begin{prop}\label{relativesequence}
Let $E$ be a profinite spectrum. Let $U\subset X$ be an open subscheme of $X$. We get a long exact sequence of cohomology groups 
$$
\cdots \to E^n_{\et}(X) \stackrel{j^{\ast}}{\to} E^n_{\et}(U) \stackrel{\partial}{\to}  E^n_{\et}(X,U) \stackrel{i^{\ast}}{\to} E^{n+1}_{\et}(X) \stackrel{j^{\ast}}{\to} E^{n+1}_{\et}(U) \to \cdots
$$
where $j:\hEt U \hookrightarrow \hEt X$ and $i:(\hEt X,\oslash) \hookrightarrow (\hEt X, \hEt U)$ denote the natural induced inclusions.
\end{prop}
\begin{proof}
Since $\hEt U \hookrightarrow \hEt X \to \hEt(X)/\hEt(U)$ is isomorphic to a cofiber sequence in $\hShp$, this is just the usual long exact sequence of $\Hom$-groups induced by a cofiber sequence in a simplicial model category.
\end{proof}
\begin{prop}\label{hominv}
Let $E$ be a profinite spectrum. Let $X \in \Sm/k$. The projection $p:X\times \A^1 \to X$ induces an isomorphism $p^{\ast}:E_{\et}^{\ast}(X) \stackrel{\cong}{\to} E_{\et}^{\ast}(X\times \A^1)$.
\end{prop}
\begin{proof}
This is clear since $p$ induces a weak equivalence in $\hSh$, see Theorem \ref{LhEt}, and hence induces isomorphisms on cohomology theories.
\end{proof}
\begin{cor}\label{EH}
Let $E$ be a profinite spectrum. Let $V \to X$ be an $\A^n$-bundle over $X$ in $\Sm/k$. Then $p^{\ast}:E_{\et}^{\ast}(X) \stackrel{\cong}{\to} E_{\et}^{\ast}(V)$ is an isomorphism. 
\end{cor}
We may summarize these results in the following
\begin{theorem}\label{Ecohom}
Let $E$ be a profinite spectrum. The \'etale cohomology theory $E_{\et}^{\ast}(-)$ represented by $E$ satisfies the axioms of a cohomology theory on $\Sm/k$ of \cite{panin}.
\end{theorem}
\begin{proof}
We check the axioms of a cohomology theory in the sense of \cite{panin}, Definition 2.0.1.\\ 
1. Localization: This is clear from Proposition \ref{relativesequence}.\\
2. Excision: Let $e:(X',U') \to (X,U)$ be a morphism of pairs of schemes in $\Sm/k$ such that $e$ is \'etale and for $Z=X-U$, $Z'=X' -U'$ one has $e^{-1}(Z)=Z'$ and $e:Z' \to Z$ is an isomorphism. By \cite{milne} III, Proposition 1.27, we know that the morphism $e$ induces an isomorphism in \'etale cohomology $H^{\ast}(\hEt(X)/\hEt(U);\Zl) \cong H^{\ast}(\hEt(X')/\hEt(U');\Zl)$. Hence the map $\hEt(X')/\hEt(U') \to \hEt(X)/\hEt(U)$ is an isomorphism in $\hHhp$. Therefore, it induces the desired isomorphism 
$E^{\ast}(\hEt(X)/\hEt(U)) \cong E^{\ast}(\hEt(X')/\hEt(U'))$ for every \'etale cohomology theory.\\
3. Homotopy invariance: This is the content of Proposition \ref{hominv}.
\end{proof}

\begin{theorem}\label{etahss}
Let $k$ be a field and let $E$ be profinite spectrum such that each coefficient group $E^q$ is finitely generated abelian. For every scheme $X$ in $\Sm/k$, there is a convergent spectral sequence 
$$E_2^{p,q}= H_{\et}^p(X;\Zln \otimes E^q) \Longrightarrow E_{\et}^{p+q}(X;\Zln).$$
\end{theorem} 
\begin{proof}
By the hypothesis on $E^q$ the groups $\Zln \otimes E^q$ are finite. The spectral sequence above is hence the one of Proposition \ref{ahss} together with the isomorphisms  
$H^{\ast}_{\et}(X;\Zln \otimes E^q) \cong H^{\ast}(\hEt X;\Zln \otimes E^q)$
of Remark \ref{remarketalecohom}. The condition for strong convergence in Proposition \ref{ahss} is satisfied since $H^p_{\et}(X;\Zln \otimes E^q)$ is a finite group for all $p$ by \cite{milne} VI, Corollary 5.5.
\end{proof}
\subsection{Examples}
\paragraph{Profinite \'etale K-theory}
\begin{defn}\label{profinetaleKU}
We define the {\rm profinite \'etale K-theory} of a smooth scheme $X$ to be the cohomology theory represented by the profinitely completed spectrum $\hKU$. 
\end{defn}
As a first application, we consider a comparison statement for profinite \'etale K-theory and for Friedlander's \'etale K-theory in \cite{etaleK}. Let $\{X_s\} \in \pro-\Sh$ be a pro-object in $\Sh$ and let $BU$ be the simplicial set representing complex K-theory. Friedlander defines the K-theory of $\{X_s\}$ for $\epsilon=0,1$ and $k>0$ by
$$K^{\epsilon}(\{X_s\};\Zln)=\Hom_{\pro-\Hh}(\{\Sigma^{\epsilon}X_s\},\{P^nBU\wedge C(\ell^{\nu})\})$$
where $\{P^nBU\}$ denotes the Postnikov tower of $BU$ considered as a pro-object in $\Sh$ and $C(\ell^{\nu})$ is the cofibre of the multiplication by $\ell^{\nu}$ map $S^1 \to S^1$.
For a locally noetherian scheme $X$, Friedlander defines the \'etale K-theory by
$K_{\et}^{\epsilon}(X;\Zln):=K^{\epsilon}(\Et X;\Zln)$. By the Atiyah-Hirzebruch spectral sequences for $\hat{KU}$ and $K_{\et}$ of Theorem \ref{etahss} and of \cite{etaleK} respectively the following result is an immediate consequence.
\begin{prop}\label{competalektheory}
For every scheme $X$ of finite type over a separably closed field $k$ the
\'etale $K$-theory groups $K_{\et}^{\ast}(X;\Zln)$ of \cite{etaleK} are isomorphic to the profinite \'etale
$K$-theory groups $\hKU_{\et}^{\ast}(X;\Zln)$ defined above. 
\end{prop} 
\paragraph{Profinite \'etale Morava K-Theory}
\begin{defn}
We define the {\rm profinite \'etale Morava $K$-theory} of a smooth scheme $X$ to be the cohomology theory represented by the profinitely completed Morava $K$-theory spectrum $\hat{K}(n)$ for $n \geq 0$.
\end{defn}
This theory satisfies in particular the conditions of Theorem \ref{etahss}. For $n\geq 1$, the coefficients are given by $\hat{K}(n)^{\ast}=\Zl[v_n,v_n^{-1}]$, where $\ddeg~ v_n=-2(\ell^n-1)$. In particular, for $n \geq 1$, the groups $\hat{K}(n)^q$ are finite for each $q$. Hence there is a converging spectral sequence from \'etale cohomology to \'etale Morava $K$-theory for $n \geq 1$:
$$E_2^{p,q}= H_{\et}^p(X;\hat{K}(n)^q) \Longrightarrow \hat{K(n)}_{\et}^{p+q}(X).$$
For schemes over $\C$, it agrees with the usual topological Morava $K$-theory. This follows as for the \'etale cobordism below from the generalized Riemann Existence Theorem of \cite{artinmazur}.
\paragraph{Profinite \'etale cobordism} 
\begin{defn}\label{profinetaleMU}
We define the {\rm profinite \'etale cobordism} of a smooth scheme $X \in \Sm/k$, to be the \'etale cohomology theory represented by the profinite cobordism spectrum $\hMU$.
\end{defn}
\begin{prop}\label{coeffofhMUet}
Let $R$ be a strict local henselian ring. Then $$\hMU^{\ast}_{\et}(\Spec R)\cong MU^{\ast}\otimes_{\Z} \Z_{\ell}~\mbox{and}~
\hMU^{\ast}_{\et}(\Spec R;\Zln)\cong MU^{\ast}\otimes_{\Z} \Zln.$$
\end{prop}
\begin{proof}
This may be easily deduced from Example \ref{piofMU}, since $\hEt \Spec R$ is contractible. The last assertion follows since $MU^{\ast}$ has no torsion.
\end{proof}
\begin{prop}\label{hEtkelltorsion} 
The (reduced) \'etale cobordism of a finite field $k$, $\car k \neq \ell$, is given by the isomorphism $\hMU^n_{\et}(k)= \hMU^{n-1}$ and similarly for $\hMU^{\ast}_{\et}(k;\Zln)$.
\end{prop}
\begin{proof}
The profinite space $\hEt k$ has the homotopy type of the circle $S^1$ in $\hSh$.
\end{proof}

The generalized Riemann Existence Theorem of \cite{artinmazur}, Theorem 12.9, or the comparison theorem for \'etale cohomology and the Atiyah-Hirzebruch spectral sequences for complex and \'etale cobordism imply the following
\begin{prop}\label{complexvar}
Let $X$ be an algebraic variety over $\C$. Let $X(\C)$ be the topological space of complex points. For every $\nu$, there is an isomorphism $$\hMU_{\et}^{\ast}(X;\Zln)\cong MU^{\ast}(X(\C);\Zln).$$
\end{prop}
There are several other obvious applications of the Riemann Existence Theorem. For example, we can deduce a proper base change theorem or an isomorphism for separably closed field extensions etc. from the corresponding theorems for \'etale cohomology. Furthermore, we can calculate the \'etale cobordism of a smooth projective curve over a separably closed field.\\
Finally, we remark that the morphism of profinite spectra $\hMU \to H\Zln$ induced by the orientation  yields a unique map of profinite \'etale cohomology theories
$\hMU^{\ast}_{\et}(X;\Zln) \to H_{\et}^{\ast}(X;\Zln)$ for every $X$ in $\Sm/k$.

\section{Algebraic versus \'etale cobordism}
The main application of the \'etale realization functor is existence of a canonical map from algebraic to \'etale cobordism. We have to show that \'etale cobordism is an oriented theory on $\Sm/k$.
\subsection{Etale cobordism is an oriented cohomology theory}
We prove that $\hMU_{\et}^{2\ast}(-;\Zl)$ is an oriented cohomology theory in the sense of \cite{lm} on the category $\Sm/k$ of smooth quasi-projective schemes over a suitable base field. The key point is the projective bundle formula.\\
Let $k$ be a field and let $G_k$ be its absolute Galois-group. We have to restrict our attention to the fields with finite $\ell$-cohomological dimension satisfying the following property:
\begin{defn}\label{defn-no-ell-torsion}
1. A profinite space $X$ is said to be {\em without $\ell$-torsion} if the canonical map 
$H^{\ast}(X;\Zln) \to H^{\ast}(X;\Zl)$ is surjective for all integers $\nu \geq 0$.\\
2. A field $k$ is said to be {\em without $\ell$-torsion} if the canonical map in Galois-cohomology 
$H^{\ast}(G_k;\Zln) \to H^{\ast}(G_k;\Zl)$ is surjective for all integers $\nu \geq 0$.
\end{defn}
\begin{example}\label{hEtPnelltorsion}
1. Every separably closed field is without $\ell$-torsion.\\
2. A finite field has no $\ell$-torsion.\\
3. Local fields are without $\ell$-torsion if $\ell$ is different from the residue characteristic of $k$.\\
4. Let $\Pro^n_k$ be the projective space of dimension $n$ over $k$. If $k$ is a field without $\ell$-torsion then the profinite space $\hEt \Pro_k^n$ has no $\ell$-torsion in $\hSh$ for every $n$. This follows from the projective bundle formula for \'etale cohomology.
\end{example}
In the rest of this chapter we will always assume that $k$ is a field without $\ell$-torsion and with finite $\ell$-cohomological dimension. 
\begin{prop}\label{hMUofPn}
The \'etale cobordism of  $\Pro^n_k$ is given by the direct sum
$$\hMU_{\et}^{\ast}(\Pro^n_k;\Zl)\cong \oplus_{i=0}^{n-1} \hMU_{\et}^{\ast}(k;\Zl).$$ 
\end{prop}
\begin{proof}
Since $\hEt \Pro_k^n$ is a space without $\ell$-torsion, the corresponding Atiyah-Hirzebruch spectral sequence collapses at the $E_2$-term and gives an isomorphism 
\begin{equation}\label{isom2.1.9}
\hMU^{\ast}(\hEt \Pro_k^n;\Zl) \stackrel{\cong}{\longrightarrow} \oplus_s H^s(\hEt \Pro_k^n;\Zl \otimes MU^{\ast -s}),
\end{equation}
as described in \cite{dehon}, Proposition 2.1.9. The assertion now follows from the projective bundle formula for \'etale cohomology.
\end{proof}

I am grateful to Francois-Xavier Dehon for an explanation of the arguments proving the following 
\begin{prop}\label{relativeKunnethforPn}
For the projective $n$-space and every smooth scheme $X$ over $k$ the canonical map
$$\hMU^{\ast}_{\et}(\Pro^n_k;\Zl)\otimes_{\hMU^{\ast}_{\et}(k;\Zl)} \hMU^{\ast}_{\et}(X;\Zl)
\stackrel{\cong}{\longrightarrow} \hMU^{\ast}_{\et}(\Pro^n_k \times_k X;\Zl)$$
is an isomorphism.
\end{prop}
\begin{proof}
This is a classical argument. In the profinite setting, this has been worked out by Dehon. The case with $\Zl$-coefficients can be done in a similar way. It is deduced from the K\"unneth formula for the $\Zl$-cohomology of $\hEt \Pro^n_k \wedge_{\hEt k} \hEt X$. The point is that $\hMU_{\et}^{\ast}(\Pro^n_k;\Zl)$ (resp. $H_{\et}^{\ast}(\Pro^n_k;\Zl)$) is a free $\hMU_{\et}^{\ast}(k;\Zl)$-module (resp. $H_{\et}^{\ast}(k;\Zl)$-module).
\end{proof}

Let $\Oh(1) \to \Pro^n_k$ be the canonical quotient line bundle. We conclude from the proof of Proposition \ref{hMUofPn} that the orientation map $\hMU^{\ast}(\hEt \Pro^n_k) \to H^{\ast}(\hEt \Pro^n_k;\Zl)$ factors through the isomorphism (\ref{isom2.1.9}). In particular, this implies that the element $\xi_H =c_1(\Oh(1)) \in  H^2(\hEt \Pro^n_k;\Zl)$ induces an element $\xi_{\hMU} \in \hMU^2(\hEt \Pro^n_k;\Zl)$, which corresponds to a morphism in $\hSHh$:
$$\xi_{\hMU}: \hat{\Sigma}^{\infty}(\hEt \Pro^n_k) \longrightarrow \hMU/\ell \wedge S^2.$$
Now let $E \to X$ be a vector bundle over $X$ in $\Sm/k$ and let $\Oh(1)$ be the canonical quotient line bundle over $\Pro(E)$. This bundle determines a morphism 
$\Pro(E) \to \Pro^N_k$ for some sufficiently large $N$. Together with the morphism $\xi_{\hMU}$ we get an element $\xi_{\hMU} \in \hMU^2_{\et}(\Pro(E);\Zl)$.
\begin{theorem}\label{PBF}{\rm Projective Bundle Formula}\\
Let $E \to X$ be a rank $n$ vector bundle over $X$ in $\Sm/k$. Then $\hMU_{\et}^{\ast}(\Pro(E);\Zl)$ is a free $\hMU_{\et}^{\ast}(X;\Zl)$-module with basis $(1, \xi, \xi^2, \ldots, \xi^{n-1})$.
\end{theorem}
\begin{proof}
We prove the assertion first for the case of a trivial bundle on $X$. As in Lemma \ref{Kunneth}, one proves that the canonical morphism of profinite spaces $\hEt (X) \times_{\hEt k} \hEt (\Pro^n_k) \to \hEt (X\times_k \Pro^n_k)$ is a weak equivalence in $\hSh$. Hence this case follows from Propositions \ref{hMUofPn} and \ref{relativeKunnethforPn}. For the general case, since $E$ is locally trivial for the Zariski topology on $X$, it suffices to show that the theorem holds for $X$ if it holds for open subsets $X_0$, $X_1$ and $X_0 \cap X_1$, with $X=X_0 \cup X_1$. This may be checked by a standard argument using the Mayer-Vietoris-sequence.  
\end{proof}

We use Grothendieck's idea to introduce higher Chern classes for vector bundles. 
\begin{defn}\label{chernclasses}
Let $E$ be a vector bundle of rank $n$ on $X$ in $\Sm/k$ and let $\xi_{\hMU}$ be as above. Then we define the {\rm $i$th Chern class of $E$} to be the unique element $c_i(E) \in \hMU_{\et}^{2i}(X;\Zl)$ such that $c_0(E)=1$, $c_i(E)=0$ for $i>n$ and $\Sigma_{i=0}^n (-1)^i c_i(E)\xi^{n-i}=0$. 
\end{defn}
The results of Panin in \cite{panin} then imply the following
\begin{theorem}\label{transfer}
The \'etale cobordism $ \hMU_{\et}^{2\ast}(-;\Zl)$ is an oriented ring cohomology theory on $\Sm/k$ in the sense of {\rm \cite{panin}}, Definition 3.1.1. This implies that for every projective morphism $f:Y \to X$ of codimension $d$ there are natural transfer maps 
$$f_{\ast}:\hMU_{\et}^{\ast}(Y;\Zl) \to \hMU_{\et}^{\ast+2d}(X;\Zl)$$
which satisfy the axioms of a push-forward in an oriented cohomology theory in the sense of {\rm \cite{lm}}. 
\end{theorem}
\begin{remark}\label{orientation}
1. If $k$ is separably closed, there are several simplifications. First of all, the above statements all hold without $\Zl$-coefficients, i.e. the Projective Bundle and the K\"unneth Formula hold for $\hMU$ and $\hMU/\ell^{\nu}$ as well.\\
Secondly, the choice of an orientation is canonical. Using the isomorphism $\LhEt MGL \cong \hMU$, we choose $x_{\hMU} \in \hMU^2(\hEt \Pro_k^{\infty})$ as the image of the orientation $x_{MGL}:\Pro^{\infty} \to MGL \wedge \Pro^1$ under $\hEt$.\\
2. A priori, the definition of $E(\Pro^{\infty}):=\lim_n E(\Pro^n_k)$ for a cohomology theory $E$ in \cite{panin} 1.1 might differ from our definition. But since $\hMU_{\et}^i(\Pro^n_k;\Zl)$ is a finite group for every $i$ and $n$ by Proposition \ref{hMUofPn}, the $\lim^1_n$ of these groups vanishes. Hence there is a canonical isomorphism  $\hMU_{\et}^i(\Pro^{\infty}_k;\Zl)\cong \lim_n\hMU_{\et}^i(\Pro^n_k;\Zl)$ for $\Pro^{\infty}_k=\colim_n \Pro^n_k$. Hence, in the case of \'etale cobordism with $\Zl$-coefficients, Panin's and our definitions are compatbile.
\end{remark}
\subsection{Comparison with $\Omega^{\ast}$}
We consider the algebraic cobordism theory $\Omega^{\ast}(-)$ of \cite{lm}. As a corollary of Theorem \ref{transfer}, using the universality of $\Omega^{\ast}$ we get
\begin{theorem}\label{OmegatohMU}
There is a canonical morphism of oriented cohomology theories
$$\theta:\Omega^{\ast}(X) \to \Omega^{\ast}(X;\Zl) \to \hMU_{\et}^{2\ast}(X;\Zl)$$  
defined by sending a generator $[f:Y\to X] \in \Omega^{\ast}(X)$, $f:Y \to X$ a projective morphism between smooth schemes, to the element $f_{\ast}(1_Y) \in \hMU_{\et}^{2\ast}(X;\Zl)$.
\end{theorem}
If $k$ is a field of characteristic zero, Theorem 4 of \cite{lm} states that $\Omega^{\ast}(k)$ is isomorphic to the Lazard ring $\Lee^{\ast}$. It is conjectured that this isomorphism holds for every field, see the conjecture of \cite{lm} \S 4.3.2.  
\begin{prop}\label{compofcoeff}
For every separably closed field $k$, the morphism 
$$\theta:\Omega^{\ast}(k;\Zln) \to \hMU_{\et}^{2\ast}(k;\Zln)$$ 
is surjective.\\
If we suppose in addition that the conjecture of \cite{lm} is true for $k$, then $\theta$ is an isomorphism. 
\end{prop}
\begin{proof}
Consider the canonical map $\Phi:\Lee^{\ast} \otimes \Zln \to \Omega^{\ast}(k;\Zln)$. It is split injective for every field, see \cite{lm}, Corollary 4.3.3. When we compose this map with $\theta: \Omega^{\ast}(k;\Zln) \to \hMU_{\et}^{2\ast}(k;\Zln)$, we get the canonical map 
$\Lee^{\ast} \otimes \Zln \stackrel{\Phi}{\to} \Omega^{\ast}(k;\Zln) \stackrel{\theta}{\to} \hMU_{\et}^{2\ast}(k;\Zln)$. Since this map is unique, it must be the canonical isomorphism. Hence $\theta: \Omega^{\ast}(k;\Zln) \to \hMU_{\et}^{2\ast}(k;\Zln)$ is surjective. The last assertion follows in the same way if we assume the conjecture to be true.
\end{proof}
\subsection{Comparison with $MGL^{\ast,\ast}$}
Let $k$ be a field as above. Let $V$ be a vector bundle of rank $d$ over $X$ in $\Sm/k$. We recall that the Thom space $\Thom (V) \in \LU$ of $V$ is defined to be the quotient $\Thom(V) =V/(V-i(X))$, where $i:X \to V$ denotes the zero section of $V$. We reformulate a lemma from $\A^1$-homotopy theory.
\begin{prop}
Let $V$ be a vector bundle over $X$ and $\Pro(V) \to \Pro(V \oplus \Oh)$ be the closed embedding at infinity. Then the canonical morphism of pointed sheaves $\Pro(V\oplus \Oh)/\Pro(V) \to \Thom(V)$ induces a weak equivalence in $\hShp$ via $\hEt$.  
\end{prop} 
\begin{proof}
This is the same proof as for Proposition 3.2.17 of \cite{mv} where we use the fact that $\hEt$ preserves $\A^1$-weak equivalences between smooth schemes by Theorem \ref{LhEt} and commutes with quotients. 
\end{proof}

Now we define the Thom class of $V$ in $\hMU_{\et}^{2d}(\Thom(V);\Zl)$. From the isomorphism $\Pro(V\oplus \Oh)/\Pro(V) \cong \Thom(V)$ we deduce an exact sequence induced by the cofiber sequence
$$\hMU_{\et}^{\ast}(\Thom(V);\Zl) \to \hMU_{\et}^{\ast}(\Pro(V\oplus \Oh);\Zl) \to \hMU_{\et}^{\ast}(\Pro(V);\Zl).$$
Using the projective bundle formula of Theorem \ref{PBF} this sequence is isomorphic to the exact sequence 
$$\hMU_{\et}^{\ast}(\Thom(V);\Zl) \to \hMU_{\et}^{\ast}(X;\Zl)[1,u,\ldots,u^d] \to \hMU_{\et}^{\ast}(X;\Zl)[1,\ldots,u^{d-1}].$$  
The element $x:=u^d -c_1(V)u^{d-1}+\ldots+(-1)^dc_d(V) \in \hMU_{\et}^{\ast}(X;\Zl)[1,u,\ldots,u^d]$ is sent to $u^d -c_1(V)u^{d-1}+\ldots+(-1)^dc_d(V) \in \hMU_{\et}^{\ast}(X;\Zl)[1,u,\ldots,u^{d-1}]$ which is $0$ by the definition of Chern classes. By exactness, there is an element $\thomet(V) \in \hMU_{\et}^{2d}(\Thom(V);\Zl)$ that is sent to $x \in \hMU_{\et}^{\ast}(X;\Zl)[1,u,\ldots,u^d]$. We call $\thomet(V)$ the {\em \'etale Thom class of $V$}. It corresponds to a morphism in $\hSHh$
$$\thomet(V):\hat{\Sigma}^{\infty}(\hEt \Thom(V)) \to \hMU/\ell \wedge S^{2d}.$$
We apply this argument to the tautological $n$-bundle $\gamma_n$ over the infinite Grassmannian. Since $MGL_n =\Thom(\gamma_n)$, we get a morphism in $\hSHh_2$   
$$\phi:\LhEt(MGL) \longrightarrow \hMU/\ell.$$
Recall the definition of algebraic cobordism for a scheme $X$ in $\Sm/k$
$$MGL^{p,q}(X):=\Hom_{\SHh(k)}(\Sigma^{\infty}_{\Pro^1}(X),MGL[p-2q]\wedge (\Pro^1)^{\wedge q}),$$
where $MGL[n]:=MGL\wedge (S^1_s)^{\wedge n}$, $S^1_s$ denoting the simplicial circle. One sees that the elements in degree $p,q$ are sent to elements in degree $p$ via $\phi$. We summarize this discussion in the following
\begin{theorem}\label{MGLcomparison}
There is a natural map for every $X$ in $\Sm/k$
$$\phi: MGL^{p,q}(X;\Zl) \to \hMU_{\et}^{p}(X;\Zl)$$
which defines a morphism of oriented cohomology theories.
\end{theorem}
By the universality of $\Omega^{\ast}(-)$, we get the following
\begin{cor}\label{triangle}
For every $ \in \Sm/k$, there is a canonical commutative diagram of morphisms of oriented cohomology theories 
$$\begin{array}{ccc}
\Omega^{\ast}(X;\Zl) & \stackrel{\theta_{MGL}}{\longrightarrow} &  MGL^{2\ast,\ast}(X;\Zl)\\
\theta_{\hMU} \searrow &  & \swarrow \phi \\
 & \hMU_{\et}^{2\ast}(X;\Zl). & 
\end{array}$$
\end{cor}
\subsection{The Galois action on \'etale cobordism}
Let $k$ be a field of characteristic $p \neq \ell$ and let $\ok$ be a separable closure of $k$. If $X$ is a scheme over $k$, there is a natural action on $\Xok=X \otimes_k \ok$ of the Galois group $G_k:=\Gal(\ok/k)$ of $k$. By naturality $G_k$ also acts on $\hEt \Xok$.\\ 
Now $\Omega^{\ast}(\ok)$ is generated by classes $[\pi:Y \to \Spec \ok]$ of smooth projective schemes $Y$ over $\ok$, see \cite{lm} \S 2.5.4. Hence in order to determine the action of the absolute Galois group on \'etale cobordism, it suffices to know the action on each such $Y$. Since $\ok$ is separably closed, the induced morphism $\sigma: Y \to Y^{\sigma}$ induces an isomorphism on \'etale cohomology with $\Zl$-coefficients, where $Y^{\sigma}$ is the twisted $\ok$-scheme with structure morphism $\sigma \circ \pi$. Hence $\sigma$ induces an isomorphism of homotopy classes $\sigma:\hEt Y_{\ok} \stackrel{\cong}{\to} \hEt Y^{\sigma}$ in $\hHh$.  
This implies the following 
\begin{prop}\label{Galoisaction}
The action of $G_k$ on $\hMU^{\ast}_{\et}(\ok;\Zln)$ is trivial.
\end{prop}
Next, we calculate the profinite \'etale cobordism groups of a local field, i.e. either a finite extension of the field $\Q_p$ or a finite extension of the field of formal power series $\F((t)))$ over a finite field of characteristic $p$. We assume $p\neq \ell$ and we denote by $q=p^f$ the number of elements in the residue class field of $k$. Let $\ell^{\nu_0}=(q-1,\ell^{\nu})$ be the greatest common divisor of $q-1$ and $\ell^{\nu}$. 
\begin{cor}\label{hMUofQp}
Let $k$ be a local field. For all $n$, the (reduced) profinite \'etale cobordism groups with $\Zln$-coefficients of $k$ are given by 
\[
\hMU^n_{\et}(k;\Zln) = \left \{ \begin{array}{r@{\quad: \quad}l}
                                          \Z/\ell^{\nu_0}\otimes MU^{n-2} & n ~\mathrm{even} \\
                                          \Z/\ell^{\nu_0} \otimes MU^{n-1} & n ~\mathrm{odd}.
                                          \end{array} \right. 
\] 
\end{cor} 
\begin{proof}
Since $\Zln \otimes MU^t$ is a finitely generated free $\Zln$-module with trivial $G_k$-action, we may identify the Galois cohomology groups $H^i(k;\Zln \otimes MU^t)$ with $H^i(k;\Zln)\otimes MU^t$. The assertion follows from the local Tate Duality for Galois cohomology and the Atiyah-Hirzebruch spectral sequence.
\end{proof}

\subsection{Inverting the Bott element}
Let $H^p(X;\Z/n(q))$ denote the motivic cohomology of a smooth scheme $X$ over a field $k$. For $\Spec k$ there is an isomorphism $H^0(\Spec k;\Z/n(1)) \cong \mu_n(k)$ with the group of $n$-th roots of unity in $k$. Assuming that $k$ contains an $n$-th root of unity $\zeta$, we have a corresponding motivic Bott element $\beta_n \in H^0(\Spec k;\Z/n(1))$. Levine has shown in \cite{bott} that motivic $\Z/n$-cohomology of a smooth scheme over $k$ agrees with \'etale $\Z/n$-cohomology after inverting the Bott element. I am grateful to Marc Levine for an explanation of his ideas.\\ 
Furthermore, Hopkins and Morel announced the construction of a motivic Atiyah-Hirzebruch spectral sequence. It is the slice filtration spectral sequence conjectured by Voevodsky in \cite{openproblems} from motivic cohomology with coefficients in the ring $MU^{\ast}$ to algebraic cobordism
\begin{equation}\label{MGLAHSS}
E_2^{p,q,2i}=H^{p,q}(X,MU^{2i}) \Longrightarrow MGL^{p+2i,q+i}(X)
\end{equation}
with differentials $d_{2r+1}:E_{2r+1}^{p,q,2i} \to E_{2r+1}^{p+2i+1,q+i,2i-2r}$.\\
We consider the spectral sequence for $\Spec k$. Since $E_2^{p,1,q}(\Spec k)$ is concentrated in degrees $p=0$ and $p=1$, we deduce $MGL^{0,1}(k)\cong k^{\times}$. For $\Z/n$-coefficients, the exact sequence for coefficients implies that we get an isomorphism 
$MGL^{0,1}(k;\Z/n)\cong \mu_n(k)$ and, via the spectral sequence, the motivic Bott element defined above, is sent to an induced {\em Bott element} $\beta_n \in MGL^{0,1}(k;\Z/n)$.

Let us now suppose that $k$ is algebraically closed, $n=\ell^{\nu}$ and $\car k \neq \ell$. In particular, $k$ contains an $\ell^{\nu}$-th root of unity $\zeta$. It defines an element $\zeta \in H_{\et}^0(\Spec k;\mu_{\ell^{\nu}})$. The element $\zeta \cdot 1_{\hMU/\ell^{\nu}}$ induces via the Atiyah-Hirzebruch spectral sequence an element $\zeta_{\hMU/\ell^{\nu}} \in \hMU_{\et}^0(\Spec k;\Zln)$. The multiplication with $\zeta$ yields an isomorphism of spectral sequences and hence multiplication with $\zeta_{\hMU/\ell^{\nu}}$ is an isomorphism on \'etale cobordism. Furthermore, the map $\phi$ sends $\beta_{\ell^{\nu}}$ to $\zeta_{\hMU/\ell^{\nu}}$. This implies that $\phi$ induces a localized map $\phi: MGL^{\ast,\ast}(X;\Zln)[\beta ^{-1}] \to \hMU_{\et}^{\ast}(X;\Zln)$.\\
Finally, motivic cohomology may be represented in $\SHh(k)$ by the Eilenberg-Mac Lane spectrum $H_{\Zln}$, defined in degree $n$ by the group completion of the symmetric products 
$$K(\Zln(n),2n) = (\coprod_{d\geq 0} S^d((\Pro^1)^{\wedge n}))^+\otimes \Zln.$$ 
Using the simplicial version of the Dold-Thom theorem and the K\"unneth isomorphism for symmetric products in \'etale cohomology over an algebraically closed field of \cite{sga4}, Expos\'e XVII, Th\'eor\`eme 5.5.21, one shows as before that the \'etale realization of $H_{\Zln}$ is isomorphic in $\hSHh$ to the profinite Eilenberg-Mac Lane spectrum representing $\Zln$-cohomology. Hence the \'etale realization defines a map from motivic $\Zln$-cohomology to \'etale $\Zln$-cohomology. With a fixed isomorphism $\Zln(1) \cong \Zln$, then this map is the unique map of ortiented cohomology theories. Furthermore, localization at $\beta$ is exact. Consequently, the \'etale realization yields a map of Atiyah-Hirzebruch spectral sequences, whose $E_2$-terms agree by Theorem 1.1 of \cite{bott}. This implies the following

\begin{condtheorem}\label{Bottconjecture}
Let $X$ be a smooth scheme of finite type over an algebraically closed field $k$ with $\car k \neq \ell$. If we assume the existence and convergence of the above Atiyah-Hirzebruch spectral sequence $(\ref{MGLAHSS})$ from motivic cohomology to algebraic cobordism, then $\phi$ is an isomorphism
$$\phi: (\oplus_{p,q}MGL^{p,q}(X;\Zln))[\beta ^{-1}] \stackrel{\cong}{\longrightarrow} \oplus_{p,q}\hMU_{\et}^{p}(X;\Zln).$$
\end{condtheorem}

\bibliographystyle{amsplain}

Mathematisches Institut, Universit\"at M\"unster, Einsteinstr. 62, D-48149 M\"unster\\
e-mail: gquick@math.uni-muenster.de

\end{document}